\documentclass{kluwer}

 \newtheorem{Definition}{Definition}
 \newtheorem{Proposition}{Proposition}
 \newtheorem{Theorem}{Theorem}
 \newtheorem{Corollary}{Corollary}
 \newtheorem{Lemma}{Lemma}
 \newenvironment{Examples}{\vspace{0.1cm}\noindent{\bf Examples.}}{\vspace{0.1cm}}
 \newenvironment{Proof}{\vspace{0.1cm}\noindent{\bf Proof:}}{\hfill\rule{2mm}{2mm}\vspace{0.1cm}\par\noindent}
 \newenvironment{Remark}{\vspace{0.1cm}\noindent{\bf Remark.}}{\vspace{0.1cm}\par}
 \newenvironment{General Assumption}{\vspace{0.1cm}\noindent{\bf General Assumption.}}{\vspace{0.1cm}\par}

 \def\eps{\varepsilon}
 \def\vol{\mathrm{vol}}
 \def\rplus{\lbrack 0, \infty )}

 \def\Scal{\mathrm{Scal}}
 \def\tr{\mathrm{tr}}
 
 \def\W{\mathbb{W}}
 \def\N{\mathbb{N}}
 \def\Q{\mathbb{Q}}
 \def\P{\mathbb{P}}
 
 \def\R{\mathbb{R}}
 \def\FF{\mathcal{F}}
 \def\PP{\mathcal{P}}
 \def\QQ{\mathcal{Q}}
 \def\DD{\mathcal{D}}
 \def\om{\omega}
 \def\Om{\Omega}
 \def\Ric{\mathrm{Ric}}
 
 \def\G{\mathcal{G}}

\begin{document}

\begin{article}
\begin{opening}
\title{Chernoff's Theorem and Discrete Time Approximations of Brownian Motion 
 on Manifolds}
\author{Oleg G. \surname{Smolyanov}}
\institute{Faculty of Mechanics and Mathematics, Moscow
 State University, Russia}
\author{Heinrich \surname{v. Weizs\"acker}}
\institute{Fachbereich Mathematik, Technische Universit\"at
 Kaiserslautern, Germany }
\author{Olaf \surname{Wittich}}
\institute{Mathematisches Institut, Universit\"at T\"ubingen, Germany}
\begin{abstract} Let $(S(t))_{t \ge 0}$ be a one-parameter family of positive
integral operators on a locally compact space $L$. For a 
possibly non-uniform partition of $[0,1]$ define a finite 
measure on the path space $C_L[0,1]$ by using a) $S(\Delta t)$ for the 
transition between any two consecutive partition times of distance $\Delta t$ 
and b) a suitable 
continuous interpolation scheme (e.g. Brownian bridges or geodesics). 
If necessary normalize the result to get a 
probability measure. We prove a version of Chernoff's theorem of semigroup 
theory and tightness results which yield convergence in law of such measures 
as the partition gets finer. In particular let $L$ be a closed smooth 
submanifold without boundary of a manifold $M$. We prove convergence of 
Brownian motion on $M$, conditioned to visit $L$ at all partition times, to
a process on $L$ whose law has a density with respect to Brownian motion on 
$L$ which contains scalar, mean and sectional curvatures terms. 
Various approximation schemes for Brownian motion on $L$ are also given.
\end{abstract}
\end{opening}

\newpage

\tableofcontents


 \section{Introduction}\label{1}

 This paper is an extension of earlier work published in two conference 
proceedings
 (\cite{Smolyanov-Weizsaecker-Wittich00a} and 
 \cite{Smolyanov-Weizsaecker-Wittich03}.
 The classical Chernoff Theorem states roughly that in the strong sense
 $S(t/r)^r \to e^{t DS}$, where $S = (S(t))_{t \ge 0}$ is a strongly
 continuous operator family on a Banach space and $DS$ its
 derivative at $t=0$, cf. \cite{Ethier-Kurtz}, p. 32, Cor. 6.6. 
 This means, under some technical assumptions,
 that two operator families $S$, $S'$ with $DS = DS'$
 yield in the limit the same semigroup. We call such families {\em
 Chernoff equivalent}. We are particularly interested in the case
 of families of positive integral operators on a smooth closed manifold $L$ 
without
 boundary. The iterations $S(t/r)^r f$  are then given by iterated integrals 
 with finite and positive kernels. Slightly extending a result from 
 \cite{Smolyanov-Weizsaecker-Wittich00a} we actually give a version of 
 Chernoff's Theorem (Proposition \ref{Chernoff}) for
 nonuniform partitions $\PP$ of the time interval $[0,t]$.\\

 For every starting point $x\in L$, time horizon $t > 0$ and
 partition $\PP = 0 = t_0 < \cdots < t_r = t$, the finite family of operators 
 $S(\Delta t_k) \circ \cdots \circ S(\Delta t_1)$, $ 0 \le
 k \le r$ defines a finite measure $\P^x_{\PP}$ on $L^{\PP}$. It
 can then be extended to a measure on the path space $C_L[0,t]$ by
 continuous interpolation, either deterministically or using suitable
 conditional distributions, or more simply on
 the path space $D_L[0,t]$ if we extend the discrete paths as stepfunctions.\\

 This construction, described in detail in Section \ref{3}, depends
 on the family $S$ and on the interpolating measures. It will
 be called the {\em pinning construction}. We are interested in the
 possible weak limit of these measures as the mesh of the partition tends to 
 $0$. The Chernoff result implies convergence of the finite dimensional 
marginals.\\ 
 
 If the family $S$ is Chernoff equivalent to a Feller semigroup we
 prove in section \ref{Fellersection}, in extension of similar
 results in \cite{Ethier-Kurtz}, the tightness  of the resulting step 
 processes over the path space $D_L[0,t]$. If the limit process has continuous 
 paths and interpolation by geodesics is used then even tightness and hence 
convergence in law over $C_L[0,1]$ follows. 
 Based on a Large-Deviation result from \cite{Wittich03b}, we allow also 
 interpolation by Brownian bridges. If $L\subset M$ is isometrically embedded 
 into another Riemannian manifold $M$, we may even interpolate by Brownian
 bridges in the ambient manifold. (Theorem \ref{mainpinning}).\\

 A particular feature of our results is the fact that the positive
 operators $S_t$ need not to be normalized in the sense that the
 associated measures $q(t,x,-)$ are allowed to have finite total mass 
$\not= 1$. For example if their densities are the restrictions of
 probability densities on the larger manifold $M$ to $L$, then the
 measures $\P^x_{\PP}$ constructed above are not probability
 measures.\\ 

 There are two different normalization procedures. First
 one can pass at the beginning from $S_t$ to the associated
 probability kernels by normalizing each $q(t,x,-)$. This gives a
 family $\tilde{S}$ which may be Chernoff equivalent to a Markov
 semigroup. In section 5 we give a couple of examples which are all
 equivalent to the heat semigroup on $L$. In order to verify this
 equivalence we use a detailed study of the short time behaviour of
 Gaussian integrals from \cite{Smolyanov-Weizsaecker-Wittich03}
 which we review in section 4. After these preparations the remaining work 
 lies in the local differential geometry which one needs for the Taylor 
 expansions of the normalization coefficients.\\

 The second possibility is to renormalize the measures $\P^x_{\PP}$
 {\it after} their construction. This is the content of section 6.
 The corresponding operator
 families are no longer Chernoff equivalent to the above Markov semigroup but
 nevertheless the resulting limit measures are equivalent to the law of
 the Markov process obtained by the first normalization procedure.
 In all situations studied in section 5 one can even calculate the 
 Radon-Nikodym density with respect to Wiener measure over $C_L[0,1]$ 
 explicitely by a combination of curvature terms. 
 This includes as a special case a new proof of a result of 
 \cite{Andersson-Driver}. However for us the most important 
 example is the Brownian motion on a larger manifold $M$ which is conditioned 
 to visit the embedded manifold $L$ at all partition times. In this case we get
 the following main result of this paper.

 \begin{Theorem}\label{surfacemeasure} Let $x\in L$ and $L(\eps):= \lbrace
 x\in M\,:\,d_M(x,L)<\eps\rbrace$ the {\em tubular
 $\eps$-neighborhood} of $L$ and $\PP_k$ be a sequence of
 partitions of the unit interval with mesh $\vert \PP_k\vert\to 0$
 as $k\to \infty$. Then the limit law of the conditional Brownian
 motions on $M$
 \begin{equation}\label{defmu}
     \lim_{k\to\infty} \lim_{\eps\to 0} \W^x_M (d\om \, \vert \, \om(t_i)\in L(\eps)\, , \, t_i
     \in\PP_k) = \mu_L^x
 \end{equation}
 exists and is equivalent to the
 Wiener measure $\W_L^x$ on the submanifold with density
 \begin{equation}\label{Dichte}
     \frac{d\mu_L^x}{d\W_L^x}(\om) = c \cdot e^{\int_0^1\left(\frac{1}{4}
 \Scal_L - \frac{1}{8} \vert\tau_{\phi}\vert^2 -
 \frac{1}{12}(\overline{R}_{M/L} +
 \overline{\Ric}_{M/L}+\Scal_M)\right)(\om(s))ds}
 \end{equation}
 where the constant $c$ normalizes to total mass $1$. \end{Theorem}

 {\bf Remark.} This limit law can be interpreted as the infinite dimensional 
surface measure which is induced by $\W^x_M$ on the submanifold $C_L[0,1]$ of 
the path space $C_M[0,1]$. This interpretation is underlined by the fact that 
the two limits in (\ref{defmu}) may be interchanged without affecting the 
result. The proof of 
this fact requires however very different techniques and is even more involved.
It has been established in \cite{Sidorova-Smolyanov-Weizsaecker-Wittich04} 
for $M = \R^m$ with tools from stochastic analysis. In this case 
\cite{Sidorova04} provides an explanation of this intechangeability of the two
limiting procedures which is independent of our present methods. The 
general case will be treated with perturbation theoretical methods in 
\cite{Wittich05} and \cite{Sidorova-Smolyanov-Weizsaecker-Wittich05}.

 \section{Chernoff's Theorem}\label{2}

 In this section we formulate a slightly extended form of the Chernoff type 
 result from \cite{Smolyanov-Weizsaecker-Wittich00a}. In contrast to the usual 
 versions we consider convergence along nonuniform partitions of the time 
 parameter. In \cite{Smolyanov-Weizsaecker-Wittich00a} we treated only 
 contractions whereas in the result below more
 general bounded operators are allowed.

 \begin{Definition}\label{properdef} Let $B(V)$
 denote the space of bounded linear operators on the Banach space $V$.
 A strongly continuous family $S : \rplus \rightarrow B(V)$ with
 $S(0)=1$ is called {\bf proper}, if
 \begin{equation}\label{normbound}
 \Vert S(t) \Vert = 1 + O(t)
 \end{equation}
 as $t \downarrow 0$, and if there is an operator $(A,\DD(A))$
 which is the generator of a strongly continuous semigroup $(e^{tA})$ on $V$ such 
 that
 \begin{equation}\label{chass}
 \frac{S(t) - I}{t}f \rightarrow Af
 \end{equation}
 as $t \downarrow 0$, for all $f \in V$ of the form $f = e^{aA}g$ with $a > 0$ 
 and $g \in V$. The operator $(A,\DD(A))$ will be also denoted by $DS$.
 \end{Definition}

 {\bf Remarks} 1. Clearly the operator $A$ is uniquely determined by the
 family $S(t)$. Therefore the notation $DS$ is justifed.
 However many proper families may lead to the same operator
 $A$. This will be discussed in the next subsection.\\
 2. The vectors of the form $e^{aA}g$ form a core for $A$, i.e. they 
 are dense in the domain $\DD(A)$ with respect to the graph norm. 
 We do not require the convergence (\ref{chass}) for {\it all} $f \in \DD(A)$ 
 since we do not need it in the following proof and our formally weaker 
 assumption is easier to verify in the situations to be studied
 later. However in contrast to the usual form of Chernoff's Theorem
 we do not know whether it suffices for the following result to
 require the above convergence for $f \in \DD$ where $\DD$ is an
 {\it arbitrary} core of the operator $A$. Note that the operators
 $(S(t) - I)/t$ may fail to be uniformly bounded in the graph norm
 of $\DD(A)$. Note also in the usual Chernoff expression
 $S(t/r)^r$ the factors commute whereas in general in a product of
 the form $S(t_1) \cdots S(t_r)$ they do not. Therefore it is not
 surprising that in comparison with the usual statement we need
 slightly stronger assumptions.

 \begin{Proposition}\label{Chernoff} Let $S(t),\ t \ge 0$ be a proper
 family of linear operators on $V$ and let $A = DS$.
 Let $(t_i^ n) _{1\le i \le r_n},\ n \in \mathbb{N},$ be
 positive numbers such that $\sum^ {r_n}_{i=1} t_i^n
 \longrightarrow t \ge 0$ and $\max_i t_i^ n \longrightarrow 0$.
 Then, for all $f \in V$, we have
 $$ S (t_{1}^ n) \ldots S(t_{r_n}^n)\ f \rightarrow e^{tA} f $$
 as $n \rightarrow \infty$.
 \end{Proposition}

 \begin{Proof} In the case of contractions this is Proposition 3 of 
 \cite{Smolyanov-Weizsaecker-Wittich00a}. The same proof works under the weaker
 assumption (\ref{normbound}) since this condition implies the existence of a 
 number $q > 0$ such that $\Vert S(t)\Vert \le e^{qt}$ for all sufficiently 
 small $t$ and thus
 \begin{equation}\label{expnormbound}
 \Vert S(t_1) \cdots S(t_r) \Vert \le e^{q \sum_{i=1}^r t_i}
 \end{equation}
 as soon as $\max_i t_i $ is small enough.
 \end{Proof}

 The usual application of Chernoff's Theorem is the construction of
 discrete {\em semigroups} $U^{(n)}$, $n=1,2,...$ approximating a 
 strongly continuous semigroup $e^{tA}$ in the sense of \cite{Kato80}, 
 (3.9), p. 511. The key observation that the limit semigroup only depends on 
 the derivative of the family $S$ at $t=0$, stays valid in our setting. 
 This will be formalized in the next section.

 \subsection{Chernoff Equivalence}\label{2.1}

 Let $\Pi$ denote the set of all proper families. We focus on the map which 
 assigns to each proper family its corresponding contraction semigroup. The
 map $P : \Pi\rightarrow \Pi$ given by $$P(S)_t := e^{t DS}$$ we call 
 {\it Chernoff map}. $P$ maps proper families onto the subset 
 $\Sigma\subset\Pi$ of strongly continuous contraction semigroups. 
 $\Sigma$ remains pointwise fixed under $P$. We are interested in the 
 attracting domains for each fixpoint, i.e. the set of preimages of a given 
 semigroup.

 \begin{Definition}\label{ChEq}
 Two proper families $S,T\in\Pi$ are called Chernoff equivalent if
 one of the two equivalent conditions holds
 \begin{itemize}
   \item[(i)]  $P(S)=P(T)$,
   \item[(ii)] $DS=DT$.
 \end{itemize}
 \end{Definition}

 The following simple criteria for Chernoff equivalence will be
 applied in the sequel.

 \begin{Lemma}\label{suff}
 Let $S$ with $DS=(A,\DD{A})$. Let $T = (T(t))$ be a family of
 operators satisfying the bound (\ref{normbound}) and
 \begin{equation}\label{oclose}
    \Vert T(t) f - S(t) f \Vert = o (t)
 \end{equation}
 for $t \downarrow 0$ and all $f$ of the form $f = e^{aA}g,\ a > 0,\ g \in V$, then $T$ is also proper and Chernoff equivalent to $S$.
 \end{Lemma}

 \begin{Proof}
 Let $f = e^{aA}g,\ g \in V$ be given. Then by our assumption (\ref{oclose}) and
 by (\ref{chass}) we get
 $$
 \lim_{t\to 0} \frac{T(t) f - f}{t} = \lim_{t\to 0}\frac{S(t) f -
 f}{t} = Af.
 $$
 Thus by definition \ref{properdef} $DT = (A,\DD(A))$, i.e. $T$ is Chernoff
 equivalent to $S$.
 \end{Proof}

 \begin{Lemma}\label{suff2}  Let $S$ be a proper family.
 Let $c = (c(t))$ be a family of operators satisfying $\Vert c(t) -
 I \Vert = o(t)$. If $T(t) = c(t) S(t)$ then $T$ and $S$ are Chernoff
 equivalent.
 \end{Lemma}

 \begin{Proof} This follows from Lemma \ref{suff} and the estimate
 $$\Vert T(t) - S(t) \Vert = \Vert (c(t) - I)S(t) \Vert = o(t)\Vert S(t) \Vert = o(t).$$
 \end{Proof}

 \section{The Pinning Construction}\label{3}

 We want to use the results of the preceding section to construct resp. approximate
 laws of Markov processes on path spaces.  The term pinning construction in the title of this section is motivated by the particular example to be studied later,
 namely by the problem of pinning a Brownian motion on a manifold $M$ down to
 a submanifold $L \subset M$ by restricting the heat kernel
 of $M$ to $L$ which gives a proper family of integral operators on $L$.\\

 The construction of the approximating processes is done in two steps: Given a
 partition $\PP$ of the time interval and a family $S$ of positive
 integral operators on $L$ we construct in natural way measures
 $\P^x_{L,\mathcal{P}}$ on the finite product $L^{\PP}$. These can be extended to
 measures $\P^x_{\PP}$ on the space of $M$-valued functions on $[0,1]$ where
 $M \supset L$ by allowing the path to make excursions into the surrounding space
 $M$ in the partition intervals.
 In this second step there is a choice of various interpolation
 schemes, ranging from step functions over geodesic interpolation to interpolation
 via Brownian bridges.\\

 Having established the approximation of the limit semigroup in the functional
 analytic sense by Chernoff's theorem (which in probabilistic language amounts
 to convergence of the finite dimensional marginals) the additional
 problem is to prove tightness of the laws of the approximating processes. This
 again has two parts. The first part is tightness under stepwise or otherwise
 trivial interpolation. In the next subsection we improve quite general 
 criteria of \cite{Ethier-Kurtz} for tightness
 over the space $D_M[0,1]$ which extend similar results for $C_{\R^m}[0,1]$ 
 given by \cite{Stroock-Varadhan79}. 
 The second part is to prove tightness over $C_M[0,1]$ for the interpolation 
by Brownian bridges which requires additional tools in the case of general 
manifolds.

  {\bf Notation.} Denote by $\vert \PP \vert$ the 
{\em mesh} of a partition $\PP$, i.e. the length of the
 longest partition interval. For locally compact separable metric space $L$ let
 $\widehat{C}(L)$ denote the Banach space of continuous 
 functions vanishing at infinity. Note that a Feller semigroup on
 $\widehat{C}(L)$ is the transition semigroup of a 
 (strong Markov) process with c\`adl\`ag paths, cf. \cite{Ethier-Kurtz}, 
 Theorem 3.2.7, p. 169.

 \subsection{Discrete Time Approximations of Feller processes}
\label{Fellersection}

 Theorem \ref{Feller} below shows that if a 
proper family of integral operators on a locally compact separable metric 
space is Chernoff equivalent to a Feller semigroup then,
 for any sequence of partitions whose mesh converges to $0$, the associated 
 measures over the piecewise constant paths converge weakly for the Skorokhod 
topology to the law of the corresponding Feller process. A similar result for 
uniform partitions and families of Markov operators is Theorem 3.2.6 in 
\cite{Ethier-Kurtz}. 

\begin{Definition}\label{Pdef} Let $(L,d)$ be a metric space and let $S$ a 
one-parameter family of integral operators
$$     
S(t) f (x) = \int \rho(t,x,dy) f(y)
$$
where  $\rho(t,x,-)$ forms a finite nonnegative Borel measure on $L$ for all $t > 0$, $x\in L$ such that $\rho(t,x,L)$ is bounded in $x$ for each $t$. 
We call $S$ a {\em pinning family}.\\
Let $ \PP := \lbrace 0= t_0<t_1\cdots<t_r = 1 \rbrace$ be a
partition of the interval $[0,1]$. For every $x \in L$ we
define a finite measure on the discrete time 'path space' $L^{\PP}$ by
\begin{eqnarray}
&& \P^x_{L,\mathcal{P}}(A_1\times\cdots\times A_r)
\label{defpmass}\\ \nonumber &=& \int_{A_1} \cdots \int_{A_r}
\rho(t_1 - t_0,x,d y_1) \cdots \rho(t_r - t_{r-1}, y_{r-1}, d y_r ).
\end{eqnarray}
Also we denote by $\overline{\P}^x_{L,\mathcal{P}}$ the unique measure
on the space $D_L[0,1]$ of c\`adl\`ag $L$-valued paths whose projection to
$L^{\PP}$ is $\P^x_{L,\mathcal{P}}$ and which is concentrated on the set of
paths which are constant on each of the partition intervals $[t_i,t_{i+1})$.
\end{Definition}

 In other words $\overline{\P}^x_{L,\mathcal{P}}$ is the image of
 $\P^x_{L,\mathcal{P}}$ under the canonical embedding of $L^{\PP}$ into
 $D_L[0,1]$.\\

 {\bf Remarks.} 1. The measures $\rho(t,x,-)$ and hence the measures 
 $\P^x_{L,\mathcal{P}}$ and $\overline{\P}^x_{L,\mathcal{P}}$ are only finite 
 but not necessarily probability measures.  Nevertheless we shall use the 
 usual topology of weak convergence of measures which is induced by the 
 duality with bounded continuous functions.
 If the $\rho(t,-,-)$ actually are probability kernels
 we can interprete $\P^x_{L,\mathcal{P}}$ and $\overline{\P}^x_{L,\mathcal{P}}$
 as the laws of two Markov process starting in $x\in L$ with state space $L$ 
 and time parameter set ${\PP}$ and $[0,1]$, respectively.\\
 2. Up to section 5 the families $S$ will be proper. In section 6 we apply 
Definition \ref{Pdef} to non proper families.\\

 \begin{Theorem}\label{Feller}
 Let $(L,d)$ be a locally compact and separable metric space. 
 Let $S = (S(t))$ be a pinning family. Suppose that $S$
 is proper and Chernoff equivalent on $\widehat{C}(L)$ to the Feller semigroup 
 $e^{tA}$. Then for every $x \in L$ and for every sequence
 of partitions $\PP_k$ with $|\PP_k| \rightarrow 0$ the associated measures
 $\overline{\P}^x_{L,\mathcal{P}_k}$
 converge weakly over the space $D_L[0,1]$ to the law of the Feller process $X$
 starting in $x$ with generator $A$.
 \end{Theorem}

 A key observation is the following routine connection between Markov chains
 and martingales.

 \begin{Lemma}\label{martingal} Let
 $Y_i, i = 0,\cdots,r$ be a (non homogeneous) Markov chain with transition 
 operators $S_i$. Then for every real bounded measurable function 
 $f$ the process
 \begin{equation}\label{martdis} f(Y_j) - \sum_{i <j}(S_i - I)f(Y_i)
 \end{equation}
 $j = 0,\cdots, r$ is a martingale with respect to the natural filtration of $(Y_i)$. 
 \end{Lemma}

 \begin{Proof} (of the theorem) 1. First let us reduce the proof to the 
normalized case where each $S(t)$ is a Markov transition operators. 
Since $e^{tA}$ is a
 Markovian semigroup we have $1 = e^{tA}1$ and $A\ 1 = 0$.
 Thus by the assumption of Chernoff equivalence and (\ref{chass}) we get
 $$\Vert S(t)1 - 1 \Vert = o(t).$$
 The  total mass $\rho(t,x,L)$ of the measure $\rho(t,x,-)$ is equal to $S(t)1 (x)$. Hence if we consider the probability measure
 $$\tilde{\rho}(t,x,-) = \frac{\rho(t,x,-)}{\rho(t,x,L)}$$
 the associated Markov operator $\tilde{S}(t)$ differs in operator norm from
 $S(t)$ only by the order $o(t)$. In particular according to Lemma \ref{suff}
 the families $S$ and $\tilde{S}$ are Chernoff equivalent. From (\ref{defpmass}) we see that
 replacing the family $S$ by the family $\tilde{S}$ changes the associated
 measure $\P^x_{L,\PP}$ in total
 variation norm only by $o(|\PP_k|)$. This implies a fortiori that the original
 sequence and the corresponding sequence of probability measures have the same
 weak limit.\\
 2. Now assume that the measures under consideration are probability measures.
 For each $k$ we consider the family of operators
 $T_k(s)$ which is defined by
 \begin{equation}\label{Tk} T_k(s) = \frac{S(t_{j+1}-t_j) - I}{t_{j+1}-t_j}\ for\  s\in [t_j, t_{j+1}).
 \end{equation}
 Let $X_{\PP_k}^x$ denote the process with law $\overline{\P}^x_{L,\mathcal{P}_k}$.
 Then
 \begin{eqnarray}\nonumber
 \sum_{i < j}(S(t_{i+1} - t_i) - I)f(X_{\PP_k}^x(t_i)) &=& \sum_{i < j}T_k(t_i)f(X_{\PP_k}^x(t_i))(t_{i+1} - t_i)\\
 &=& \int_0^{t_j} T_k(s)f(X_{\PP_k}^x(s)) ds.\label{sumint}
 \end{eqnarray}
 For every function of the form $f = e^{aA}g$ we know from the assumption
 of Proposition \ref{Chernoff} that $T_k(s)f(z) \rightarrow Af(z)$ uniformly in
 $s\in [0,1], z\in L$ as $k \rightarrow \infty$. Fix $k \in \N$ and $a > 0$, 
 and let $f$ be of the form $f = e^{aA}g,$ $ g \in \widehat{C}(L)$. Put 
  $$M^f_k(t) = f(X^x_{\PP_k}(t_j)) - \sum_{i <j}(S(t_{i+1} - t_i) - I)
    f\big(X^x_{\PP_k}(t_i)\big)\,\,\, \mathrm{for}\,\ t\in [t_j,t_{j+1})$$
and 
 \begin{equation}\label{martdec}
 Z^f_k(t)  = M^f_k(t) -  f(X^x_{\PP_k}(t)).
 \end{equation}
Then by  Lemma \ref{martingal} the process $M^f_k(t)$ is a martingale with 
respect to the natural filtration $(\mathcal{G}^k_t)$ of $(X^x_{\PP_k}(t))$,
because $\mathcal{G}^k_t = \mathcal{G}^k_{t_j}$ for $t\in [t_j,t_{j+1})$. 
Moreover by (\ref{sumint}) and the above uniform convergence there is a 
finite deterministic constant $C_f$ which depends on $f$ such that 
 $$\sup_k |Z^f_k(t)| = \sup_k |\int_0^{t_j} T_k(s)f(X_{\PP_k}^x(s)) ds| 
\le C_f t.$$
 The set of all functions of the form $f = e^{aA}g$ is uniformly dense in 
 $\widehat{C}(L)$. Thus Theorem 9.4 of chapter 3 of \cite{Ethier-Kurtz} can be 
applied to the algebra $C_a = \widehat{C}(L)$ and we conclude that for each 
$f \in \widehat{C}(L)$ the
 sequence of processes $f \circ X^x_{{\PP}_k}$ is uniformly tight in
 $D_{\R}[0,1]$.\\
 Moreover the finite dimensional marginals of the processes $X^x_{\PP_k}$
 converge to the the corresponding marginals of the process $X^x$ by
 Proposition \ref{Chernoff}. Hence Corollary 9.3 of Chapter 3 of \cite{Ethier-Kurtz} gives the assertion.
 \end{Proof}

 In our applications we are interested in weak convergence over the space
 of continuous rather than c\`adl\`ag paths. For the corresponding transfer 
 between these settings the following Lemma is useful:

 \begin{Lemma}\label{margstraff} Let $(L,d)$ be a separable metric space.
 For a function $ \varphi : [0,1] \rightarrow L$ and $\delta > 0$ let
 $$ w(\varphi, \delta) := \sup\{ d (\varphi(s), \varphi(t)) : s,t \in [0,1], 
|s-t| < \delta \}.$$
 Let the sequence $(\P_k)$ of finite measures over $D_L[0,1]$ 
converge weakly to a finite measure $\P_0$ which is concentrated on 
$C_L[0,1]$. Then for all $ \varepsilon > 0 $ one has, as $\delta \downarrow 0$,
 \begin{equation}\label{fluccont}
 \limsup_{k \in N} {\P_k} \{ \omega
 \in D_L[0,1] : w(\omega, \delta) > \varepsilon
 \} \longrightarrow 0.
 \end{equation}
 \end{Lemma}

 \begin{Proof} Since the total mass of the $\P_k$ converges to the total mass
of $\P_0$ we may assume that we deal with probability measures. A sequence in 
$D_L[0,1]$ which converges in the topology of 
this space to an element of $C_L[0,1]$ actually converges uniformly on $[0,1]$,
 cf. \cite{Ethier-Kurtz}, Lemma 3.10.1. By Skorokhod representation 
\cite{Ethier-Kurtz}, Theorem 3.1.8,
 there are a probability space $(\Om, \FF, \Q)$ and processes
 $X^k, k \ge 0$ on this space such that for each $k$, $X^k$ has law $\P_k$ and 
 $X^k \rightarrow X^0$ a.s. in $D_L[0,1]$ and hence also a.s. uniformly on
 $[0,1]$.
 Since $X^0$ has a.s. (uniformly) continuous paths we have for each 
$\varepsilon > 0$
 $${\Q}\{w(X^0, \delta) > \varepsilon
 \} \longrightarrow 0.$$
 Because of the a.s. uniform convergence of the paths this implies 
  $$\limsup_{k \in N}{\Q}\{w(X^k, \delta) > \varepsilon
 \} \longrightarrow 0.$$
 This is equivalent to (\ref{fluccont}).
 \end{Proof}

 \subsection{Continuous Interpolations}

 We assume now that $L$ is embedded into another locally compact space $M$. 
 We construct a net of measures
 $\P^x_{\mathcal{P}}$ on the path space $C_M[0,1]$ associated to
 the family $S$. It is indexed by the {\em finite}
 partitions $\mathcal{P}$ of $[0,1]$ and depends on a starting
 point $x\in L$. The marginals of these
 measures on $M^{\PP}$ are concentrated on $L^{\PP}$ and given by the measures
 $\P^x_{L,\PP}$ introduced in Definition \ref{Pdef}. In the partition intervals  we use an 'interpolation family':

 \begin{Definition}\label{Qdef}
 A family $\mathcal{Q}:=\lbrace\mathbb{Q}_{s,t}^{x,y} \,:\, x,y\in L, 0\leq s<t \leq 1
  \rbrace$ of probability measures on the path space $C_M \lbrack s,t \rbrack$
 such that
 $$\mathbb{Q}_{s,t}^{x,y}(\lbrace\omega \,:\, \omega(s)=x, \omega(t) = y\rbrace  ) = 1$$ is called {\bf interpolating family}.
 \end{Definition}

 We combine the interpolating family
 $\mathbb{Q}_{t,s}^{x,y}$ with a measure of the form $\P^x_{L,\mathcal{P}}$ 
 to arrive at a path measure as follows: Every $\om\in C_M[0,1]$ 
 can be identified with a unique $m$-tuple
 $$
 \underline{\om}:=(\om_1,\ldots,\om_r)\in C_M[0,t_1]
 \times\cdots\times C_M[t_{r-1},1]
 $$ which satisfies $\om_j(t_j)= \om_{j+1}(t_j)$ for all
 $j\in\{1,\ldots,r-1\}$. Using this identification we define the
 measure $\P^x_{\PP}$ on $ C_M[0,1] $ by
 \begin{equation}\label{PxPdef}
  \P^x_{\PP}(d\om)  = \int_{M^{\PP}} \P^x_{L,\mathcal{P}}(dy_1\times\cdots
  \times dy_r)\,\Q_{0,t_1}^{x,y_1} (d\om_1)
 \cdots \Q_{t_{r-1},1}^{y_{r-1},y_r } (d\om_r) .
 \end{equation}
 Here $ \P^x_{L,\mathcal{P}}$ is considered as a measure on $M^{\PP}$. 
Evaluation at the partition points gives a canonical projection
 $\pi_{\PP}$ of the path space $C_M[0,1]$ to $M^{\PP}$. Then
 under this projection, the measure $\P^x_{\PP}$
 has the marginal measure $\P^x_{L,\mathcal{P}}$. In particular both measures 
have the same total mass.\\

We give three examples for possible interpolating families.

 \begin{Examples} {\bf (1) L-Geodesic Interpolation:} Let $L$ be a Riemannian 
manifold. The family  $\mathcal{Q}$ is
 given by the point mass $\Q^{x,y}_{s,t} := \delta_{\gamma^L_{x,y,s,t}}$, where
 $$\gamma^L_{x,y,s,t}(u):= \gamma^L_{x,y}(\frac{u-s}{t-s}d_L(x,y))$$ 
 and $\gamma^L_{x,y}$ is an
 arbitrary shortest geodesic in $L$ connecting $x$ an $y$
 parametrized by arc length. The measures constructed by
 $L$-geodesic interpolation are supported by the path space $C_L[0,1]$.

 \noindent{\bf (2)  M-Geodesic Interpolation:} Assume that $L$ is
 isometrically embedded into the manifold $M$. The family  $\mathcal{Q}$ 
is given by the point mass $\Q^{x,y}_{s,t} := \delta_{\gamma^M_{x,y,s,t}}$, 
where 
$$\gamma^M_{x,y,s,t}(u):= \gamma^M_{x,y}\big(\frac{u-s}{t-s}d_M(x,y)\big)$$ 
and $\gamma^M_{x,y}$ is an arbitrary shortest
 geodesic in $M$ connecting $x$ and $y$ parametrized by arc length.

 \noindent{\bf (3) Brownian Bridge Interpolation:} Here
 $\Q^{x,y}_{s,t}$ denotes the measure of a {\em Brownian bridge} in
 $M$ starting at time $s$ at $x$ and ending up at time $t$ at $y$.
 \end{Examples}

 For the sake of definiteness we fix now 
 the general setting for the sequel.
 Nevertheless many of our general arguments could be adapted to other similar 
 situations. 

 \begin{General Assumption} We assume that the Riemannian 
manifold $L$ is smooth, closed (i.e. compact without boundary) of 
dimension $l$ and isometrically embedded in to the Riemannian manifold $M$ of 
dimension $m$.
\end{General Assumption}
 
 Here is the consequence of Theorem \ref{Feller} in our context. 
 This will be the central general tool in the last two sections.  

 \begin{Theorem}\label{mainpinning} Under the 'General Assumption' above let
 $S$ be a proper pinning family which is Chernoff equivalent
on the Banach space $C(L)$ to the semigroup $(e^{tA})$ of a diffusion 
processes, i.e.
a (strong Markov) Feller process with continuous paths on $L$. Let $\QQ$ be 
either $L$-geodesic, $M$-geodesic or Brownian
 bridge interpolation or another interpolating family for which the
 assumptions of Lemma \ref{BLdist} below hold. Then for 
 every $x \in L$ and for every sequence
 of partitions $\PP_k$ with $|\PP_k| \rightarrow 0$ the associated measures
 $\P^x_{{\PP}_k}$ obtained by (\ref{PxPdef})
 converge weakly over the space $C_M[0,1]$ to the law of the process $X$
 starting in $x$ with generator $A$.
 \end{Theorem}

 \begin{Proof} (1) We begin with $M$-geodesic interpolation. Fix $\delta > 0$. 
 Choose $k$ large enough such that $\vert \PP_k \vert < \delta$. 
 Let $s < t$ with $t-s < \delta$ and a path $\omega$ be given which is 
 geodesic in the intervals of the partition 
$\PP_k = \{t_0 < \cdots < t_{r_k}\}$.
 Choose the indices $l,u$ such that 
 $t_l := \max\lbrace \tau \in\PP_k : \tau \leq s \rbrace$ and
 $t_u := \min\lbrace \tau \in\mathcal{P} : t \leq \tau \rbrace$. By
 the construction of geodesic interpolation, we have
 \begin{eqnarray}\label{geodest}
     & & d_M (\om (s),\om (t)) \\ 
\nonumber &\leq & d_M(\om (s), \om (t_{l+1})) + 
d_M(\om (t_{l+1}), \om (t_{u-1})) + d_M(\om (t_{u-1}, \om (t)) \\
\nonumber &\leq & d_M(\om (t_{l-1}), \om (t_l)) + d_M(\om (t_{l+1}), 
\om (t_{u-1})) + d_M( \om (t_{u-1}), \om (t_{u}))\\
\nonumber &\leq& 3 \max \{d_M(\om(t_i),\om(t_j)) : 
t_i, t_j \in \PP_k, |t_j - t_i| < \delta\}.
 \end{eqnarray} 
The measures $\overline{\P}^x_{L,\mathcal{P}_k}$ of Definition \ref{Pdef} 
and $\P^x_{\PP_k}$ defined in (\ref{PxPdef}) have the same marginal measure 
$\P^x_{L,\PP_k}$ on $M^{\PP_k}$ resp $L^{\PP_k}$. Let us denote by $w_M$ and $w_L$ respectively the modulus of continuity introduced in Lemma \ref{margstraff} 
with respect to the metrics $d_M$ resp $d_L$. Then 
\begin{eqnarray*} &&\P^x_{\PP_k} \{\om \in C_M[0,1] : 
w_M(\om,\delta) > \eps\}\\
&\leq & \P^x_{L,\PP_k} \{\om \in L^{\PP_k} : \max \{d_M(\om(t_i),\om(t_j)) : 
|t_j - t_i| < \delta\} > \eps/3 \} \\
&\leq & \P^x_{L,\PP_k} \{\om \in L^{\PP_k} : \max \{d_L(\om(t_i),\om(t_j)) : 
|t_j - t_i| < \delta\} > \eps/3 \} \\
&\leq & \overline{\P}^x_{L,\PP_k} \{\om \in D_L[0,1] : w_L(\om,\delta) > \eps/3 
\}. 
\end{eqnarray*} 
According to Theorem \ref{Feller} the sequence 
$(\overline{\P}^x_{L,\mathcal{P}_k})_k$ converges weakly over the space 
$D_L[0,1]$ to the law of the continuous process $X$. Thus by Lemma 
\ref{margstraff} we get for each $\eps > 0$
\begin{equation}\label{straffkrit}
\lim_{\delta \to 0} \limsup_{k \to \infty} 
\P^x_{\PP_k} \{\om \in C_M[0,1] : w_M(\om,\delta) > \eps\} = 0.
\end{equation}
This implies tightness of the sequence $(\P^x_{\PP_k})_k$ over the space 
$C_M[0,1]$ and clearly the law of $X$ is the only possible limit point.\\
(2) The proof for $L$-geodesic interpolation is completely analogous taking $d_L$, $w_L$ instead of $d_M$, $w_M$.\\
 (3) For the Brownian bridge interpolation the result now follows from  
 Lemma \ref{BLdist} and Proposition \ref{BrBri} below.
 \end{Proof}
 
  We prove that Brownian bridge interpolation leads to the same 
limit measure as geodesic interpolation. If $M=\R^m$ a much easier proof of 
this fact was given in \cite{Smolyanov-Weizsaecker-Wittich00a}.

\begin{Lemma}\label{BLdist} Let  $\mathcal{Q}$ be an interpolating family 
   such that for $\eps > 0$
 small enough, we have for all $\alpha > 0$ an increasing function
 $g_{\alpha} : \lbrack 0,\alpha\rbrack \to \R^+_0$ with $\lim_{u\to
 0} g_{\alpha}(u) = 0$ such that
 $$
 \Q_{s,t}^{x,y}(\Gamma^M_{x,y,s,t} (\alpha))\leq (t-s) g_{\alpha}
 (t-s)
 $$
 uniformly in $$R_{\eps}:=\lbrace(x,y)\in M\times M\,:\, d_M
 (x,y)<\eps\rbrace$$ where $$\Gamma^M_{x,y,s,t} (\alpha) := \lbrace
 \omega\in C_M\lbrack s,t\rbrack\,:\,\sup_{s\leq u\leq t} d_M
 (\om (u),\gamma^M_{x,y,s,t}(u))\geq\alpha \rbrace.$$ Then the
 sequence $\W^x_{\PP_k}$ constructed from $S$ and $\mathcal{Q}$
 converges to the same limit measure as the sequence constructed
 from $S$ by M-geodesic interpolation.
 \end{Lemma}

 \begin{Proof} Let $\nu_k$ and $\rho_k$ be measures on $C_M[0,1]$
 constructed from $\mathcal{Q}$ and $M$-geodesic interpolation
 respectively. For each path $\omega
  \in C_M[0,1]$ let $\varphi_k(\omega)$ be a $M$-geodesic interpolation of the
 restriction $\omega_{|\mathcal{P}_k}$. Then we have the relation
 $\rho_k = \nu_k \circ \varphi_k^{-1}$. Moreover
 \begin{eqnarray*}
 && \nu_k\{ \omega: \sup_{0 \le u \le 1}d_M\big(\omega(u),\varphi_k(\omega)(u)\big) > \alpha\}\\
 &\le& \sum_{l=1}^{r_k} \nu_k\{ \omega :\sup_{t_{l-1} \le u \le t_l}d_M\big(\omega(u),\varphi_k(\omega)(u)\big) > \alpha\}\\
 &\le& \sum_{l=1}^{r_k} \int \mathbb{P}^x_{L,\mathcal{P}_k}(dy)
 \Q^{y_{l-1},y_l}_{t_{l-1},t_l}
 \big\{\Gamma^M_{y_{l-1},y_l,t_{l-1},t_l}(\alpha)\big\}\\
 &\le&
 \overline{\mathbb{P}}^x_{L,\mathcal{P}_k}\{w(\omega,|\mathcal{P}_k|)
 > \varepsilon\} +  \sum_{l=1}^{r_k} \sup_{d_M(y_{l-1},y_l) \le
 \varepsilon} \Q^{y_{l-1},y_l}_{t_{l-1},t_l}
 \big\{\Gamma^M_{y_{l-1},y_l,t_{l-1},t_l}(\alpha)\big\}\\
 &\le&  \overline{\mathbb{P}}^x_{L,\mathcal{P}_k}\{w(\omega,|\mathcal{P}_k|)
 > \varepsilon\} +  \sum_{l=1}^{r_k} (t_l - t_{l-1})g_\alpha(t_l - t_{l-1}).\\
 \end{eqnarray*}
 Here both terms are arbitrarily small for large $k$ and small
 $\alpha$: The first by Lemma (\ref{margstraff}) and the second due
 to our assumption. Hence for every uniformly continuous function
 $f:C_M[0,1] \rightarrow \R$ we get
 \[ \lim_{k\rightarrow \infty} | \int f(\omega)\ d\nu_k -\int f(\omega)\ d\rho_k|
 \le \int |f(\omega) - f(\varphi_k(\omega))| \ d\nu_k = 0\] which
 implies that the sequences $(\nu_k)$ and $(\rho_k)$ have the same
 limit points in law.
 \end{Proof}

 \begin{Remark} 1. Instead of M-geodesic interpolation we could
 have used in Lemma \ref{BLdist} any other interpolating family as reference 
for which the pinning measures are known to be tight. In the proof of Theorem 
\ref{mainpinning} the geodesic interpolation was used to get the
second inequality in (\ref{geodest}). Clearly a uniform estimate of 
the form
$$ d_M (\om (s),\om (t)) \leq h\big(d_M(\om(t_i),\om(t_{i+1}))\big) \ for\ s,t \in 
(t_i,t_{i+1})$$
where $\lim_{\eps \to 0} h(\eps) = 0$ would have been sufficient. It is easy to
construct other interpolation schemes on more general metric spaces which 
satisfy such a condition.\\
2. By LeCam's Theorem, see \cite{Dudley}, 11.5.3 Theorem, p. 316, 
convergence of the respective sequence implies its uniform tightness.
 \end{Remark}

 Finally, we use a uniform Large-Deviation result about Brownian
 bridges \cite{Wittich03b}
 to conclude the corresponding part of Theorem \ref{mainpinning}
 from Lemma \ref{BLdist}.

 \begin{Proposition}\label{BrBri} The Brownian bridge interpolation family
 $\mathcal{Q}$ on $M$ satisfies the assumption of Lemma \ref{BLdist} with
 $g_{\alpha} (u) := 2 e^{-\chi \alpha^2/u}/u$ for some $\chi >0$.
 \end{Proposition}

 \begin{Proof} As in the proof at the end of the previous section, let
 $r_M := \inf_{x\in L} r_M (x)$ where $r_M(x)>0$ is the largest
 number such that the geodesic balls $B(x,r)$ are {\em strongly
 convex} in the sense of \cite{doCarmo}, 3.4, p. 74 for all $r < r_M
 (x)$. Let $\eps < r_M/2$. By the Large Deviation result from
 \cite{Wittich03b}, there are $\chi
 > 0$ and $\alpha_0
 > 0$ such that for all $0 < \alpha < \alpha_0$ we have
 $$
 \Q_{s,t}^{x,y}(\Gamma^M_{x,y,s,t} (\alpha))\leq 2 \exp \left(-
 \frac{\chi\,\alpha^2}{t-s}\right)
 $$
 as $t-s\to 0$, uniformly in $R_\varepsilon$. This completes the proof.
 \end{Proof}

 \section{Gaussian Integrals}\label{Gauss}

 In the sequel, proper families as described in section \ref{2}
 will be constructed by families of integral operators. In order to
 compute the derivative at zero -- and therefore the Chernoff
 equivalence class -- of such a family we first review some facts
 about the short time asymptotic of Gaussian integrals from 
 \cite{Smolyanov-Weizsaecker-Wittich03}. We introduce
 a degree $d$ on the space of space-time polynomials such that the
 short-time contribution of a monomial $p$ either vanishes or is of
 order $t^{d(p)}$. Using this notion, we reformulate Wick's formula
 in an algebraic way and conclude Corollary \ref{bedg}. It states
 that in our situation the only relevant terms are of homogeneous
 degree one. For the proofs of the results in this section cf. 
 \cite{Smolyanov-Weizsaecker-Wittich03}, p. 351-354. 

 \subsection{Wick's Formula}\label{Wick}

 Let $t>0$. By Fubini's theorem and using
 the fact that the Gaussian integral solves the heat equation, we
 obtain the following result also known as {\it Wick's formula}.

 \begin{Lemma}\label{wick}
 For $k\in\mathbb{Z}\times\mathbb{N}_0^n$ define 
 \begin{equation}\label{defdk}
 d(k):= k_0+\frac{1}{2}(k_1+\ldots + k_n).
 \end{equation}
 Let $$p_k (t,\xi) = t^{k_0}\xi_1^{k_1}\ldots \xi_n^{k_n}$$ be a monomial in $(t,\xi)$ 
 such that $d(k)\geq 0$. Then for
 \begin{equation}\label{defgt}
   \mathcal{G}_t(p_k) := \int_{\mathbb{R}^n} \frac{e^{-\frac{\vert \xi
   \vert^2}{2t}}}{\sqrt{2\pi t}^n} p_k(t,\xi)\ d\xi
 \end{equation}
 we obtain
 $$
   \mathcal{G}_t(p_k) =\left\lbrace \begin{array}{lc}
   0 & (k_1,\ldots,k_n) \notin (2\mathbb{N}_0)^n \\ t^{d(k)} \prod_{i=1}^n (k_i - 1)!!
   & \mathrm{else} \\
   \end{array}\right..
 $$
 Here we use the standard notation
 $(2n-1)!!=1\cdot3\cdot\ldots\cdot(2n-1).$
 \end{Lemma}

 Let now $\mathcal{L}$ denote the real algebra generated by all
 monomials of the form $p_k(t,x)$, $k\in\mathbb{Z}\times
 \mathbb{N}^{n}$. The map
 $$D:\mathcal{L}\rightarrow \frac{1}{2}\mathbb{Z}$$ given by $D(p_k):= d(k)$
 induces a {\it grading}
 $\mathcal{L}=\oplus_{s\in\frac{1}{2}\mathbb{Z}} \mathcal{L}^s,$
 where $\mathcal{L}^s :=\langle p_k : d(k)=s \rangle$ denotes the
 subspace of homogeneous elements of degree $s$. We consider the
 {\it associated filtration} by ideals
 $\mathcal{I}^r:=\oplus_{s\geq r} \mathcal{L}^s.$ Let furthermore
 $f\mapsto \lbrack f \rbrack$ denote the quotient map
 $$q:\mathcal{I}^0\rightarrow \mathcal{I}^0/\mathcal{I}^{3/2}$$ and
 $Q$ the projection onto the subalgebra generated by monomials
 $p_k$ with $(k_1,\ldots,k_n)\in(2\mathbb{N}_0)^n$.

 If in particular $p_k$ is a monomial with $k_0 = 0$ then $Q\lbrack p_k \rbrack = 0$ unless $p_k \in\lbrace 1, x_1^2,...,x_n^2\rbrace$. The 
 following Proposition about the short time asymptotic of Gaussian integrals is 
 a simple consequence of Lemma \ref{wick}.

 \begin{Proposition}\label{wickalg}
 Let $f\in \mathcal{I}^0$. Then $\lim_{t\to 0} \G_t (f)$ exists and
 we have asymptotically
 \begin{equation}
 \mathcal{G}_t (f) = \mathcal{G}_t(Q\lbrack f \rbrack) + o(t)
 \end{equation}
 as $t\to 0$. More explicitely, if $f=f_0+f_{1/2}+f_1+\ldots$ is the
 decomposition of $f$ into homogeneous elements, then
 \begin{equation}\label{relevant}
 \G_t (f) = \G_t(Qf_0) + \G_t(Qf_1) + o(t).
 \end{equation}
 \end{Proposition}

 \begin{Remark}
  $q$ is a ring homomorphism, whereas $Q$ is not.
 \end{Remark}

 In the sequel we will use the following fact concerning quotients
 of Gaussian integrals.

 \begin{Corollary}\label{bedg}
 Let $f,h\in\mathcal{L}$ have the homogeneous components
 \begin{itemize}
 \item[(i)] $f=f_0 + f_{1/2} + f_1$ and $f_0\in\mathbb{R}$ is
 constant. \item[(ii)] $h=1 + h_1$.
 \end{itemize}
 Then, as $t \downarrow 0$,
 \begin{equation}
  \frac{\G_t(fh)}{\G_t(h)} -
  \frac{\G_0(fh)}{\G_0(h)}= \G_t(f_1) + o(t).
 \end{equation}
 \end{Corollary}

 In order to make use of the above discussion we have to show how it applies 
 to more general situations. We observe first that the polynomial short time
 asymptotic of a Gaussian integral is in some sense independent of
 the domain of integration. From this we draw the following conclusion for the 
 polynomial short time asymptotic of more general functions (see \cite{Smolyanov-Weizsaecker-Wittich03}, Corollary 3):

 \begin{Corollary}\label{Tayl}
 Let $U\subset\mathbb{R}^n$ be a bounded open neighbourhood of the
 origin.  For $r\in \mathbb{N}$ and $f\in C^r (U)\cap C (\overline{U})$ denote 
 the Taylor polynomial of $f$ up to order $r-1$ (the ($r-1$)-jet) around $0$ by 
 $\hat{f}$. For every $k\in\mathbb{Z}$
 such that $t^k\hat{f}\in\mathcal{I}^0$ we have asymptotically as
 $t\to 0$
 $$
    \int_U\frac{e^{-\frac{\vert \xi \vert^2}{2t}}}{\sqrt{2\pi
    t}^n}t^k f(\xi)d\xi=\mathcal{G}_t (t^k \hat{f}) + O(t^{k+\frac{r}{2}})
 $$
 where the constant in the error term depends on $f$ only via the maximal 
 Taylor coefficient of $f$ of order $r$ in a small neighbourhood of $0$.
 \end{Corollary}
 
\begin{Remark} In the subsequent asymptotic computations of integrals, we will denote a $C^r$-function and its Taylor expansion by the same symbol omitting the hat introduced above. For example, for a function $f\in C^r(L)$ such that $t^k \hat{f}\in \mathcal{I}^0$, $\lbrack t^k f\rbrack$ denotes the equivalence class of $t^k \hat{f}$ in $\mathcal{I}^0/\mathcal{I}^{3/2}$.
\end{Remark}

 \section{Proper Families Equivalent to the Heat Semigroup on
 Manifolds}\label{properex}

 \subsection{Introduction}\label{properexintro}

 Again let $L$ be a closed smooth manifold and $\vol_L$ the corresponding Riemannian volume measure.
 Our examples in this section follow the same pattern: 
 Consider a family of smooth 
 integral kernels $q_t (x,y)\in C^{\infty}_+
 (L\times L)$, $t >0$ and the associated operators
\begin{equation}\label{V} 
S(t)f(x) = \int_L q_t (x,y) f(y)\vol_L (dy).
\end{equation}
on the Banach space $C(L)$.  We introduce the corresponding (normalized)
Markov operators
 \begin{equation}\label{U1}
     T(t) f(x) = \frac{\int_L q_t (x,y) f(y)\vol_L (dy)}{\int_L q_t (x,z) \vol_L
     (dz)}.
 \end{equation}
We also 
assume resp. verify that
there is a function $D \in C(L)$ such that the denominator in (\ref{U1}) 
satisfies 
 \begin{equation}\label{equ}
   b(t,x):=  \int_L q_t (x,y) \vol_L (dy) = e^{t D(x)} + o(t)
 \end{equation}
uniformly in $x \in L$ as $t \downarrow 0$. 

\begin{Proposition}\label{TtoB} Under the assumption (\ref{equ}) the  
family $(T(t))$ given by (\ref{U1}) is proper if and only if the operator 
family $(B(t))$ defined by 
 \begin{equation}\label{U2}
     B(t) f(x) = \int_L e^{-t D(x)}\,q_t (x,y) f(y)\vol_L (dy)
 \end{equation}
is proper and in this case they are Chernoff equivalent.
\end{Proposition} 

\begin{Proof} We have $B(t)f(x) = c(t)T(t)f(x)$ where $c(t)$ is the operator 
of multiplication with the function
 $$
 x \mapsto b(t,x) e^{- t D(x)}
 $$
 By (\ref{equ}) and Lemma \ref{Kern} below we
 can apply Lemma \ref{suff2} to the operators $c(t)$ and get the assertion.
\end{Proof}

Actually the corresponding semigroup will be always the heat semigroup 
$(e^{\frac{t}{2}\Delta_L})$.  

 \begin{Lemma}\label{Kern} Let $k\in C(L\times L)$, $k\geq 0$ and consider 
 the integral operator
 $$
     I_k f (x) := \int_L k(x,y) f(y) \vol_L(dy).
 $$
 Then $I_k : C(L) \rightarrow C(L)$ is a bounded operator with norm
 $$
     \Vert I_k \Vert = \sup_{x\in L}\left\vert\int_L k(x,y)\vol_L(dy)
     \right\vert=\Vert I_k 1 \Vert_\infty.
 $$
 \end{Lemma}

 Therefore we can infer from Theorem \ref{mainpinning} convergence of 
the mesures constructed by the pinning construction starting with the either 
of the families $(T(t))$ or $(B(t))$ to the law of the diffusion process generated by $A$.

 \subsection{First Examples}\label{firstex}
 In this subsection we follow, with a couple of minor corrections, 
 the exposition in \cite{Smolyanov-Weizsaecker-Wittich03} in order to 
 motivate the subsequent calculations in subsections \ref{pseudogauss} -
 \ref{Heatkernel}.  We fix the Laplace-Beltrami-operator to be non-positive as 
our choice of sign 
 and consider the heat semigroup $(e^{t\Delta_L/2})$ on the Banach space 
$C(L)$.\\ 

 \noindent Let $d_L (-,-)$ denote the distance function on $L$. Consider the 
 pseudo-Gaussian kernel 
 $$k_t(x,y) = \frac{1}{\sqrt{2 \pi t}^l} e^{-\frac{d_L (x,y)^2}{2t}}.$$
 The first step is to consider the family of associated (normalized) 
 Markov-operators
 \begin{equation}\label{AD1}
  T(t)f(x) :=\frac{\int_L e^{-\frac{d_L (x,y)^2}{2t}}f(y)\vol_L(dy)}{\int_L e^{
  -\frac{d_L (x,y)^2}{2t}}\vol_L(dy)}.
 \end{equation}
 By the smoothness of the heat kernel on $L$, the subspace $C^3(L)\subset
 C(L)$ contains the image $e^{\frac{a}{2}\Delta_L}C(L)$ for each $a > 0$. Thus 
 we prove Chernoff equivalence of the family above and the heat semigroup
 if we can show that
 \begin{equation}
  \lim_{t\to 0} \frac{T(t) f - f}{t} = \frac12\Delta_L f
 \end{equation}
 for all $f\in C^3(L)$. As noted above, we may restrict ourselves to
 integration over an arbitrary open neighbourhood $U(x)$ instead of
 over all of $L$. We choose $U(x)$ so small that we can use the
 exponential map to construct a {\it normal coordinate
 system} $$\exp^L_x : V(0) \rightarrow U(x).$$ Let $\xi := \exp^L_x (y)$. Then we have
 $$d_L (x,y) =\vert \xi \vert$$ in
 these coordinates. Furthermore $$\vol_L(dy)=\sqrt{\det
 g}(\xi)\,d\xi$$ where $g$ is the metric tensor. Therefore we obtain
 $$
   (T(t) f - f)(x) = \frac{\mathcal{G}_t(\sqrt{\det g}\,f)}{\mathcal{G}_t(\sqrt{\det g})}
   (0) - \frac{\mathcal{G}_0(\sqrt{\det g}\,f)}{\mathcal{G}_0(\sqrt{\det g})}(0).
 $$
 But now $\sqrt{\det g}$ is infinitely differentiable due to our assumptions
 on the manifold and $f$ is in $C^3(L)$. We may thus apply Corollary
 \ref{Tayl} with $k=0$ and in the sequel we just have to consider the Taylor
 expansion of these functions up to second order. The expansion of the metric
 tensor $g$ we quote from \cite{Roe}, (1.14) Proposition, p. 8

\begin{Lemma}  \label{metrik}
In normal coordinates the Taylor expansion of $g$ is given by
$$g_{ab}(\xi) = \delta_{ab} + \frac{1}{3} R^L_{auvb}(0) \xi^{u} \xi^v + O(\vert \xi
   \vert^3),$$
where $R^L$ denotes the curvature tensor of $L$.
\end{Lemma}
 Thus we get for the volume form, cf. e.g. \cite{Smolyanov-Weizsaecker-Wittich03}, Cor. 4:

\begin{Corollary}
In normal coordinates the Taylor expansion of $\sqrt{\det g}$ is given by
$$
\sqrt{\det g}(\xi) = 1 + \frac{1}{6} R^L_{auva}(0) \xi^{u} \xi^v + O(\vert \xi\vert^3),
$$ 
where $R^L$ denotes the curvature tensor of $L$.
\end{Corollary}

 Therefore the 2-jet of $\sqrt{\det g}$ has exactly the properties required
 for the function $h$ in Corollary \ref{bedg}. So we apply Corollary 
 \ref{Tayl} and Corollary \ref{bedg} to the 2-jets of $f\circ \exp^L_x$ and 
 $\sqrt{\det g}$ to get
 \begin{eqnarray*}
  (T(t) f - f)(x) &=& \mathcal{G}_t (f_1) (0) + O(t^{3/2}) \\
                 &=& \frac{t}{2} \Delta f(0) + O(t^{3/2}).
 \end{eqnarray*}
 But in normal coordinates, the Laplacian on $L$ coincides with
 $\Delta$, since we were assuming our Laplace-Beltrami operator
 always to be non-positive. Therefore we may write invariantly
 $$(T(t) f - f)(x) = \frac{t}{2}\Delta_L f (x) + O(t^{3/2}),$$ 
 for all functions $f\in C^3(L)$. This is a pointwise statement. But
 inspecting Corollary \ref{Tayl} above shows
 that due to compactness the remainder is $O(t^{3/2})$ {\em
 uniformly} on $L$. This finally implies

 \begin{Proposition}\label{AD1Ch}
 The family $(T(t))$ defined in (\ref{AD1}) is Chernoff equivalent to the 
heat semigroup on $L$. 
 \end{Proposition}

 Next we want to omit the denominator in (\ref{AD1}) and compensate it by a 
suitable modification of the kernel. For this we
 consider the short time asymptotic of this  denominator 
 $$b(t,x):=\int_L e^{-\frac{d_L (x,y)^2}{2t}}\vol_L(dy).$$

 \begin{Lemma} \label{normalisierung} Let $\Scal_L$ be the scalar curvature of 
$L$. Then, uniformly in $x \in L$,
 $$ b(t,x) = \sqrt{2 \pi t}^l \left( e^{-\frac{t
 \Scal_L(x)}{6}}+O(t^{3/2})\right).$$
 \end{Lemma}

 \begin{Proof}
 Using again the Taylor expansion of the volume form, we get
  \begin{eqnarray*}
  b(t,x) &=& \sqrt{2 \pi t}^l 
\left(\mathcal{G}_t (\lbrack \sqrt{\det g}\rbrack) + O(t^{3/2})\right) \\
         &=& \sqrt{2 \pi t}^l \left(1 + \frac{t}{2} \Delta \left(\frac{1}{6} 
R^L_{auva} \xi^{u} \xi^v \right)_{\xi=0}  +O(t^{3/2})\right)\\
         &=& \sqrt{2 \pi t}^l \left(1 + 
\frac{t}{6}  R^L_{auua}(0)  +O(t^{3/2})\right) \\
         &=& \sqrt{2 \pi t}^l \left(1 - \frac{t}{6}\Scal_L(x) 
+ O(t^{3/2})\right).
  \end{eqnarray*}
 \end{Proof}

Propositions \ref{TtoB} and \ref{AD1Ch} now imply 
 
\begin{Corollary}\label{B1} The family of bounded operators defined by 
 \begin{equation}\label{AD2}
  B(t) f (x) :=\frac{1}{\sqrt{2 \pi t}^l}\int_L  e^{-\frac{d_L
  (x,y)^2}{2t} + \frac{t\Scal_L(x)}{6}}f(y)\vol_L(dy).
 \end{equation}
is also Chernoff equivalent to the heat semigroup on $L$.   
\end{Corollary}

 Finally, it should be noted that there is also a symmetric version of the
 approximating kernel. By the very same arguments as in Lemma
 \ref{normalisierung} with
 $$
  \hat{b}(t,x) := \int_L e^{-\frac{d_L (x,y)^2}{2t} 
+ \frac{t(\Scal_L(x) + \Scal_L(y))}{12}}\vol_L(dy)
 $$
 instead of $b(t,x)$ we obtain as well:

 \begin{Corollary}
 The operator family defined by 
 $$
   \hat{B}(t) f (x) :=\frac{1}{\sqrt{2 \pi t}^l} \int_L  e^{-\frac{d_L
  (x,y)^2}{2t}+\frac{t\Scal_L(x)+ t \Scal_L(y)}{12}}f(y)\vol_L(dy)
 $$
 is proper and Chernoff equivalent to the heat semigroup on $L$.
 \end{Corollary}

 \subsection{The Heat Equation on a submanifold via the restriction of a 
pseudo-Gaussian kernel}\label{pseudogauss}

 We now consider the restriction to $L$
 of a pseudo-Gaussian kernel on $M$, and prove Chernoff equivalence to the
 heat semigroup on $L$. This result was
 communicated to us with a different proof already in \cite{Tokarev}.
 We will state it in the spirit of the preceding sections as follows:

 \begin{Theorem}\label{T1} Let $\phi:L\subset M$ be an isometric embedding of 
 the closed and connected smooth Riemannian manifold $L$ into the smooth 
Riemannian manifold $M$. Let $\dim (L)=l$, $\dim(M)=m$. Let $\Scal_L$ be the 
scalar curvature of $L$. Let $\tau_\phi$ denote the tension vectorfield of the
 embedding and 
 \begin{equation}\label{Rquer}
   \overline{R}_{M/L} := \sum_{a,b=1}^l \langle R^M(e_a,e_b)
 e_b, e_a \rangle
 \end{equation}
 the partial trace of the curvature tensor of $M$ over an
 arbitrary orthonormal base of $\phi_* TL$. Then the family $(B(t))$ defined by
 \begin{equation}\label{tokarev}
 B(t) f (x) :=\frac{e^{t\big(\frac{\Scal_L}{4}
-\frac{\vert\tau_{\phi}\vert^2}{8}-\frac{\overline{R}_{M/L}}{12}\big)(x)}}{\sqrt{2
 \pi t}^l}\int_L e^{-\frac{d_M (\phi(x),\phi(y))^2}{2t}}f(y)\vol_L(dy)
 \end{equation}
 is proper and Chernoff equivalent to the heat semigroup on $L$.
 \end{Theorem}
 Specializing this result to embeddings into euclidean space we obtain the
 following statement.

 \begin{Corollary} If, in Theorem \ref{T1}, $M = \R^m$ then the family
 $$
 B_t f (x) :=\frac{e^{t\big(\frac{\Scal_L}{4} 
-\frac{\vert\tau_{\phi}\vert^2}{8}\big)(x)}}{\sqrt{2 \pi t}^l}
\int_L e^{-\frac{\vert\phi(y)-\phi(x)\vert^2}{2t}}f(y)\vol_L(dy)
 $$
 is proper and Chernoff equivalent to the heat semigroup on $L$.
 \end{Corollary}

 \subsection{Proof of Theorem \ref{T1}}\label{Tokarevproof}

 Let $x,y\in L$ and $f\in C^3(L)$. We now want to determine an
 asymptotic expression for
 \begin{equation}\label{pseudoMfam}
    S(t) f (x) :=\frac{1}{\sqrt{2 \pi t}^l}\int_L e^{-\frac{d_M (\phi(y),\phi(x))^2}{2t}}f(y)
    \vol_L (dy).
 \end{equation}
 By the arguments above we can reduce the problem to purely local
 considerations on sufficiently small neighbourhoods $U_L(x)$ and
 $U_M(\phi(x))$ which are chosen such that $\phi(L)\cap U_M(\phi(x)) = \phi(U_L(x))$. To do so, we consider local normal coordinates 
$$
 \begin{array}{c}
  \exp^L_x : V(0) \rightarrow U_L(x) \\  \exp^M_{\phi(x)}:W(0)\rightarrow
  U_M(x). \\
 \end{array}
 $$ 
Local coordinates for $L$ and $M$ are denoted by $\xi=(\xi^1,\ldots ,\xi^l)$ and
 $\eta=(\eta^1,\ldots ,\eta^m)$ respectively. The local coordinate representation
 $$(\exp^M_{\phi(x)})^{-1}\circ\phi\circ\exp^L_x: V(0) \rightarrow W(0)$$ of
 $\phi$ will be denoted by the same letter, i.e. $$ \eta = \phi(\xi) =  (\phi^1(\xi),\ldots ,\phi^m (\xi)).$$ In these local coordinates we obtain
 \begin{equation}\label{Spseudokoor}
    S(t) f (x) =\frac{1}{\sqrt{2 \pi t}^l}\int_{V(0)} e^{-\frac{\vert\phi(\xi)\vert^2}{2t}}f(\xi)
    \sqrt{\det g^L}(\xi) d\xi + O(t^{3/2}).
 \end{equation}
 We want to apply Proposition \ref{wickalg}. To do so we have to
 make sure that
 $$\lbrack h(\xi) \rbrack = \left\lbrack e^{-\frac{\vert\phi(\xi)\vert^2-\vert 
\xi \vert^2}{2t}}\sqrt{\det g^L}f(\xi)\right\rbrack $$ 
really satisfies the assumptions made there. To see this we denote by
 \begin{equation}
   \mathrm{H}_{\phi} (-,-) = \nabla^L d\phi (-,-),
  \end{equation}
  the {\it Hessian} of the map $\phi$, which coincides with the {\it second fundamental form} of the
 embedding. Now the necessary input from differential geometry can be
 summarized in the following proposition.

 \begin{Proposition}\label{hessian}
  Consider points $x,y\in L\subset M$, where $L$ is considered as
  isometrically embedded by the map $\phi$ as described above. Let $p\in U(x)$
  and $U(x)\subset L$ so small that $y$ and $x$ are joined by a unique minimizing
  geodesic $\gamma^L_{xy}$ starting at $x$. We assume $\gamma^L_{xy}$ to be
  parametrized by arc-length. Then
  \begin{equation}
   \lim_{d_L(x,y)\to 0} \frac{d_M(x,y)^2 - d_L(x,y)^2}{d_L(x,y)^4}  = 
   - \frac{1}{12}\left\Vert
   \mathrm{H}_{\phi}(\dot{\gamma}^L_{xy}(0),\dot{\gamma}^L_{xy}(0))\right\Vert^2
  \end{equation}
 \end{Proposition}

 \begin{Proof} In those local coordinates defined above, we have using the
 Taylor expansion of $\phi$ around $\xi=0$:
 \begin{eqnarray*}
 & &  \frac{d_M(x,y)^2 - d_L(x,y)^2}{d_L(x,y)^4} = \frac{\vert\phi(\xi)\vert^2 -\vert \xi \vert^2}{\vert \xi \vert^4} \\
 &=&\frac{\left\vert\frac{\partial \phi^{\alpha}}
  {\partial \xi^a}(0) \xi^{a} +\frac{1}{2}\frac{\partial^2 \phi^{\alpha}}
  {\partial \xi^a \partial \xi^b}(0) \xi^{a} \xi^b + \frac{1}{6}\frac{\partial^3 \phi^{\alpha}}
  {\partial \xi^a \partial \xi^b \partial \xi^{u}}(0) \xi^{a} \xi^b \xi^{u}\right\vert^2
  -\vert \xi \vert^2}{\vert \xi \vert^4} + O(\vert \xi \vert).
 \end{eqnarray*}
 Denote the metric on $L$ by $g^L$ and the metric on $M$ by $g^M$.
\footnote{Throughout we use freely the Einstein summation convention. 
From now on greek indices will be used for the summation from $1$ to $m$ 
and latin indices for summation from $1$ to $l$} The fact
 that $\phi$ is an isometric embedding is equivalent to the local equation
 (see \cite{Jost}, p. 29 f.)
  \begin{equation}\label{isom}
   g^L_{ab}(\xi)
   = g^M_{\alpha \beta} (\eta) 
   \frac{\partial \phi^{\alpha}}{\partial \xi^a}
                 \frac{\partial \phi^{\beta}}{\partial \xi^b}.
  \end{equation}
 Since the metric tensor at the origin of a normal coordinate system is the
 flat one, we obtain using (\ref{isom})
 \begin{eqnarray*}
 g^M_{\alpha \beta} \frac{\partial \phi^{\alpha}}
  {\partial \xi^a} \frac{\partial \phi^{\beta}}{\partial \xi^b}(0) \xi^a \xi^b 
  = g^L_{ab} (0) \xi^a \xi^b = \vert \xi \vert^2.
 \end{eqnarray*}
 Partial differentiation of (\ref{isom}) yields 
\begin{eqnarray*}&& \frac{\partial
 g^L_{ab}}{\partial \xi^u} = \frac{\partial g^M_{\alpha \beta}}{\partial
 \eta^{\rho}}  \frac{\partial \phi^{\rho}}{\partial \xi^u }\frac{\partial
 \phi^{\alpha}}{\partial \xi^a} \frac{\partial \phi^{\beta}}{\partial
 \xi^b}+g^M_{\alpha \beta} \left(\frac{\partial^2 \phi^{\alpha}}{\partial
 \xi^a\partial \xi^u} \frac{\partial \phi^{\beta}}{\partial \xi^b}   + \frac{\partial \phi^{\alpha}}{\partial
 \xi^a} \frac{\partial^2 \phi^{\beta}}{\partial \xi^b\partial \xi^u}\right).
\end{eqnarray*} 
But at the origin of a normal coordinate system the partial derivative of the metric tensor coincides with
 its covariant derivative and therefore vanishes. This implies  
$$ \delta_{\alpha\beta}\left(\frac{\partial^2 \phi^{\alpha}}{\partial
 \xi^a\partial \xi^u} \frac{\partial \phi^{\beta}}{\partial \xi^b}   + \frac{\partial \phi^{\alpha}}{\partial
 \xi^a} \frac{\partial^2 \phi^{\beta}}{\partial \xi^b\partial \xi^u}\right)(0)  = 0.$$ 
This means
 \begin{eqnarray*}
 & &  \frac{d_M(x,y)^2 - d_L(x,y)^2}{d_L(x,y)^4} \\ &=&\frac{\delta_{\alpha\beta}
 \left( \frac{1}{4}\frac{\partial^2 \phi^{\alpha}}
  {\partial \xi^a \partial \xi^b} \frac{\partial^2 \phi^{\beta}}
  {\partial \xi^u \partial \xi^v} + \frac{1}{3}\frac{\partial^3 \phi^{\alpha}}
  {\partial \xi^a \partial \xi^b \partial \xi^{u}}\frac{\partial \phi^{\beta}}
  {\partial \xi^v}  \right)(0)\xi^{a} \xi^b \xi^{u} \xi^v}{\vert \xi \vert^4} + O(\vert\xi\vert)\\
  &=& -\frac{1}{12} \delta_{\alpha\beta}\frac{\partial^2 \phi^{\alpha}}
  {\partial \xi^a \partial \xi^b} \frac{\partial^2 \phi^{\beta}}
  {\partial \xi^u \partial \xi^v}(0)  \frac{\xi^a \xi^b \xi^u \xi^v}{\vert \xi \vert^4} \\
  & & + \frac{1}{3} \delta_{\alpha\beta}\left(\frac{\partial^2 \phi^{\alpha}}
  {\partial \xi^a \partial \xi^b} \frac{\partial^2 \phi^{\beta}}
  {\partial \xi^u \partial \xi^v} + \frac{\partial^3 \phi^{\alpha}}
  {\partial \xi^a \partial \xi^b \partial \xi^{u}}\frac{\partial \phi^{\beta}}
  {\partial \xi^v}\right)(0) \frac{\xi^a \xi^b \xi^u \xi^v}{\vert \xi \vert^4} + O(\vert\xi\vert). \\
 \end{eqnarray*}
 If we now differentiate (\ref{isom}) twice we obtain at the origin
 \begin{eqnarray*}
  & &\frac{\partial^2 g^L_{ab}}{\partial \xi^u\partial \xi^v}(0)
  = \frac{\partial^2 g^M_{\alpha\beta}}{\partial \eta^{\rho}\partial \eta^{\mu}}
  \frac{\partial \phi^{\rho}}{\partial \xi^u}\frac{\partial \phi^{\mu}}{\partial \xi^v}
  \frac{\partial \phi^{\alpha}}{\partial \xi^a}\frac{\partial \phi^{\beta}}{\partial \xi^b}(0)
  + \delta_{\alpha\beta}\left(\frac{\partial^3 \phi^{\alpha}}{\partial \xi^a \partial \xi^u \partial \xi^v}
   \frac{\partial \phi^{\alpha}}{\partial \xi^b}\right. \\
   && + \left.\frac{\partial^3 \phi^{\alpha}}{\partial
   \xi^b \partial \xi^u \partial \xi^v}\frac{\partial \phi^{\alpha}}{\partial \xi^a} + \frac{\partial^2 \phi^{\alpha}}{\partial \xi^a \partial \xi^u }
   \frac{\partial^2 \phi^{\alpha}}{\partial \xi^b \partial \xi^v} + \frac{\partial^2
   \phi^{\alpha}}{\partial \xi^a \partial \xi^v }
   \frac{\partial^2 \phi^{\alpha}}{\partial \xi^b \partial \xi^u}\right)(0).
 \end{eqnarray*}
 Using this and relating the partial derivatives of the metric tensor to
 curvature by using the fact (see Lemma \ref{metrik}) that in normal
 coordinates the Taylor expansion of $g$ is given by
 \begin{equation}\label{metrik2}
 g_{ab}(\xi) = \delta_{ab} + \frac{1}{3} R_{auvb}(0) \xi^{u} \xi^v 
+ O(\vert \xi \vert^3),
 \end{equation}
 where $R$ denotes the curvature tensor, we obtain
 \begin{eqnarray*}
  & & 2\delta_{\alpha\beta} \left(\frac{\partial^2 \phi^{\alpha}}
  {\partial \xi^a \partial \xi^b} \frac{\partial^2 \phi^{\beta}}
  {\partial \xi^u \partial \xi^v} + \frac{\partial^3 \phi^{\alpha}}
  {\partial \xi^a \partial \xi^b \partial \xi^{u}}\frac{\partial \phi^{\beta}}
  {\partial \xi^v}\right)(0)\, \xi^a \xi^b \xi^u \xi^v \\
  &=& \left(\frac{\partial^2 g^L_{ab}}{\partial \xi^u\partial \xi^v}-\frac{\partial^2 g^M_{\alpha\beta}}{\partial \eta^{\rho}\partial \eta^{\mu}}
  \frac{\partial \phi^{\rho}}{\partial \xi^u}\frac{\partial \phi^{\mu}}{\partial \xi^v}
  \frac{\partial \phi^{\alpha}}{\partial \xi^a}\frac{\partial \phi^{\beta}}{\partial \xi^b}\right)(0)\, \xi^a \xi^b \xi^u \xi^v \\
  &=& 2\left(R^L_{auvb} - R^M_{\alpha\rho\mu\beta}
  \frac{\partial \phi^{\rho}}{\partial \xi^u}\frac{\partial \phi^{\mu}}{\partial \xi^v}
  \frac{\partial \phi^{\alpha}}{\partial \xi^a}\frac{\partial \phi^{\beta}}{\partial \xi^b}
  \right)(0)\, \xi^a \xi^b \xi^u \xi^v = 0 
  \end{eqnarray*}
  due to the symmetries of the curvature tensor (see \cite{Jost}, (3.3.7),
  p. 129). On the other hand, the remaining term is indeed the Hessian at the 
origin of a local normal coordinate system (see \cite{Jost}, (3.3.47), p. 138) 
and since $\xi=(\exp^L_x)^{-1} (y)$ and therefore $\xi / \vert \xi \vert = \dot{\gamma}^L_{xy}(0)$ we finally obtain our statement.
 \end{Proof}

 We now obtain the following result, first derived with a different proof in \cite{Tokarev}.

 \begin{Proposition}\label{pseudo} With the notations of Theorem \ref{T1} we 
have for $f \in C^3(L)$
\begin{equation}\label{Asympseudo}S(t)f(x) = 
 e^{tD(x)}f(x)+
  \frac{t}{2}\Delta_L f (x) + O(t^{3/2}),
\end{equation}
in particular
\begin{equation}\label{pseudonorm}\frac{1}{\sqrt{2 \pi t}^l}\int_L e^{-\frac{d_M (\phi(y),\phi(x))^2}{2t}} \vol_L (dy) =  e^{tD(x)} + O(t^{3/2})
 \end{equation} 
where 
$$D(x) = -\left(\frac{\Scal_L}{4}
+\frac{\vert\tau_{\phi}\vert^2}{8}+\frac{\overline{R}_{M/L}}{12}\right)(x).$$
 \end{Proposition}

 \begin{Proof} By (\ref{Spseudokoor}) and Proposition
 \ref{hessian}
 \begin{eqnarray*}
 & &S(t)f(x)= \frac{1}{\sqrt{2 \pi t}^l}\int_{V(0)} e^{-\frac{\vert\phi(\xi)\vert^2}{2t}}f(\xi)
    \sqrt{\det g^L}(\xi) d\xi + O(t^{3/2})\\
&=& \mathcal{G}_t\left(Q(\lbrack
 e^{-\frac{\vert\phi(\xi)\vert^2-\vert \xi \vert^2}{2t}}\rbrack \, \lbrack
 \sqrt{\det g^L}\rbrack\, \lbrack f \rbrack) \right)+O(t^{3/2}) \\ &=&
 \mathcal{G}_t\left(Q(1+\frac{\delta_{\alpha\beta}\frac{\partial^2
  \phi^{\alpha}}{\partial \xi^a \partial \xi^b} \frac{\partial^2 \phi^{\beta}}
  {\partial \xi^u \partial \xi^v}(0) \xi^a \xi^b \xi^u \xi^v}{24t})
  (1 + \frac{1}{6} R_{iuvi}(0) \xi^{u} \xi^v )\lbrack 
  f \rbrack\right)\\
  & & +O(t^{3/2}),
 \end{eqnarray*}
 since by Lemma \ref{metrik} we have
 \begin{equation}\label{det}
    \sqrt{\det g}(\xi) = e^{\frac{1}{2}\tr\log g(\xi)} = 1 + \frac{1}{6} R_{iuvi}(0) \xi^{u}
    \xi^v+ O(\vert \xi \vert^3)
 \end{equation}
 and the error $O(\vert \xi \vert^3)$ leads to an error term $O(t^{3/2})$ in 
the Gaussian integral. Now
 $$\lbrack f \rbrack = f(0) + \frac{\partial f}{\partial \xi^s}(0)\xi^s +
  \frac{1}{2}\frac{\partial^2 f}
{\partial \xi^u\partial \xi^v}(0) \xi^v \xi^u.$$ 
But since the other two factor only contain monomials of even degree, all
 contributions containing $\partial f/\partial \xi^s$ are annihilated by $Q$.
 Thus
  \begin{eqnarray*}
 & & S(t)f(x)\\
   &=& \mathcal{G}_t\left(Q(1+\frac{\delta_{\alpha\beta}\frac{\partial^2
       \phi^{\alpha}}{\partial \xi^a \partial \xi^b} 
       \frac{\partial^2 \phi^{\beta}}{\partial \xi^u \partial \xi^v}(0) 
       \xi^a \xi^b \xi^u \xi^v}{24t} + \frac{1}{6} R^L_{iuvi}(0) \xi^{u} \xi^v)
       \right)f(0) \\
   & &  + \mathcal{G}_t\left(Q(\frac{1}{2}\frac{\partial^2 f}
       {\partial \xi^u\partial \xi^v}(0)\xi^v \xi^{u})\right) + O(t^{3/2}) \\
   &=& \mathcal{G}_t\left(1 + \frac{1}{6} R^L_{iuvi}(0) \xi^{u} \xi^v\right)f(0) 
       + \frac{1}{2}\frac{\partial^2 f}{\partial \xi^u\partial \xi^v}(0)
       \mathcal{G}_t\left( \xi^v \xi^{u}\right) \\       
   & & + \,\mathcal{G}_t \left( \frac{\delta_{\alpha\beta} \frac{\partial^2\phi^{\alpha}}
       {\partial \xi^a \partial \xi^b} \frac{\partial^2 \phi^{\alpha}}
       {\partial \xi^u \partial \xi^v}(0) \xi^a \xi^b \xi^u \xi^v}{24t}\right)
       f(0) + O(t^ {3/2})\\
   &=&  e^{-\frac{t R^L_{iuui}}{6}}f(0) + \frac{1}{24t}
       \mathcal{G}_t\left(\delta_{\alpha\beta}\frac{\partial^2
       \phi^{\alpha}}{\partial \xi^a \partial \xi^b} 
       \frac{\partial^2 \phi^{\beta}}
       {\partial \xi^u \partial \xi^v}(0) \xi^a \xi^b \xi^u \xi^v\right)f(0)\\
   & & +\frac{t}{2}\Delta f(0)+O(t^{3/2}).
 \end{eqnarray*}
 The second term remains to be computed. By
 \begin{equation}\label{DeltaDelta}
  \Delta\Delta \xi^a \xi^b \xi^u \xi^v = 8(\delta_{uv}\delta_{ab}+\delta_{ua}\delta_{vb}
  +\delta_{ub}\delta_{va})
 \end{equation}
 we obtain, having in mind that the {\em tension vector field} is the trace of the second fundamental form (see \cite{Jost}, (8.1.16), p.319 for the corresponding formula in local coordinates)
 \begin{eqnarray*}
  & & \frac{1}{24t} \mathcal{G}_t\left(\delta_{\alpha\beta}\frac{\partial^2
  \phi^{\alpha}}{\partial \xi^a \partial \xi^b} \frac{\partial^2 \phi^{\beta}}
  {\partial \xi^u \partial \xi^v}(0) \xi^a \xi^b \xi^u \xi^v\right) \\
  &=& \frac{t}{24}\delta_{\alpha\beta}\frac{\partial^2
  \phi^{\alpha}}{\partial \xi^a \partial \xi^b} \frac{\partial^2 \phi^{\beta}}
  {\partial \xi^u \partial \xi^v}(0)(\delta_{ab}\delta_{uv}+\delta_{av}\delta_{ub}+
  \delta_{au}\delta_{vb}) \\
  &=& \frac{t}{24}\delta_{\alpha\beta}\left(\frac{\partial^2
  \phi^{\alpha}}{\partial \xi^a \partial \xi^a} \frac{\partial^2 \phi^{\beta}}
  {\partial \xi^b \partial \xi^b} + 2 \frac{\partial^2
  \phi^{\alpha}}{\partial \xi^a \partial \xi^b} \frac{\partial^2 \phi^{\beta}}
  {\partial \xi^a \partial \xi^b}\right)(0) \\
  &=& \frac{t}{12} \delta_{\alpha\beta} \left(\frac{\partial^2
  \phi^{\alpha}}{\partial \xi^a \partial \xi^b} \frac{\partial^2 \phi^{\beta}}
  {\partial \xi^a \partial \xi^b}-\frac{\partial^2
  \phi^{\alpha}}{\partial \xi^a \partial \xi^a} \frac{\partial^2 \phi^{\beta}}
  {\partial \xi^b \partial \xi^b}\right)(0) \\
  & & +\frac{t}{8}\delta_{\alpha\beta}\frac{\partial^2
  \phi^{\alpha}}{\partial \xi^a \partial \xi^a} \frac{\partial^2 \phi^{\beta}}
  {\partial \xi^b \partial \xi^b}(0) \\
  &=& \frac{t}{12} \left(\langle\mathrm{H}_{\phi}(e_a,
  e_b), \mathrm{H}_{\phi}(e_a,
  e_b) \rangle - \langle\mathrm{H}_{\phi}(e_a,
  e_a), \mathrm{H}_{\phi}(e_b,
  e_b) \rangle \right)(x)\\
  & & + \frac{t}{8}\vert \tau(\phi) \vert^2(x) 
 \end{eqnarray*}
 where the vectors $e_{a, a=1...,l}$ form an orthonormal base of $\phi_* T_x L$.
 But by the Gauss equations (see \cite{Jost}, Thm. 3.6.2, (3.6.7), p.151) we
 have
 \begin{eqnarray*}
  & & \langle\mathrm{H}_{\phi}(e_a,e_b), \mathrm{H}_{\phi}
  (e_a,e_b) \rangle - \langle\mathrm{H}_{\phi}(e_a,
  e_a), \mathrm{H}_{\phi}(e_b,e_b) \rangle \\
  &=& \langle R^M(e_b,e_a)e_a,e_b  \rangle -
  \langle R^L (e_b,e_a)e_a,e_b \rangle.
 \end{eqnarray*}
 Summing over $a$, $b$ yields the statement.
 \end{Proof}
 \noindent Now, since 
$$e^{\frac{t}{2}\Delta_L}f(x) = f(x) +
 \frac{t}{2}\Delta_Lf(x) + O(t^{3/2})$$ 
 Proposition \ref{pseudo} shows  
 \begin{eqnarray*}\Vert B(t)f - e^{t/2 \Delta_L}f\Vert &=& \Vert e^{-tD}S(t)f - 
e^{t/2 \Delta_L}f\Vert\\
&=& \Vert (e^{-tD}-1) \frac{t}{2}\Delta_Lf\Vert 
+ O(t^{3/2}) = O(t^{3/2}) 
\end{eqnarray*}
for functions $f\in C^3 (L)$. By Lemma \ref{suff} the family $(B(t))$ is 
 Chernoff-equivalent to the heat semigroup. This completes the proof of the 
Theorem.\\ 

{\bf Remark.} According to Proposition \ref{TtoB} the normalized version 
$(T(t))$ (see (\ref{U1})) of the 
family $(S(t))$ defined in (\ref{pseudoMfam}) is also Chernoff equivalent to 
the heat semigroup.

 \subsection{Pinning the Heat Kernel}\label{Heatkernel}

 In the case where $M$ is Euclidean space the
 kernels in the integrals of section \ref{pseudogauss} are not only pseudo-Gaussian but
 truly Gaussian, i.e. they are induced by restricting the transition kernel of
 the Brownian motion in the ambient space to $L$. This is no longer true for 
 arbitrary $M$. In order to obtain an
 analogous interpretation in the general case as well, we are now
 interested in the family obtained by using the {\em heat kernel}
 of the ambient manifold instead of the pseudo-Gaussian one. Let thus $p_t^M (x,y)$ be
 the heat kernel on $M$ and consider on $C(L)$ the operator family
 \begin{equation}\label{heat}
     P_t f (x) := \sqrt{2 \pi t}^{m-l}
\int_L p_t^M (x,y) f(y) \vol_L (dy).
 \end{equation}
 It turns out that the calculation of an asymptotic formula for (\ref{heat})
 is an application of the previous result as long as we take for granted the
 {\em Minakshisundaram-Pleijel expansion} \cite{Minakshisundaram-Pleijel} of 
the heat kernel:

 \begin{Theorem} Let $p_t^M (0,\xi)$ be the heat kernel on $M$ in a normal coordinate
 neighbourhood around $x\in M$. Then we have the following asymptotic
 expansion
 \begin{equation}
     \sqrt{2\pi t}^m\, e^{\frac{d_M
     (0,\xi)^2}{2t}}p_t^M (0,\xi) = \left.\det\right.^{-1/4} g^M (\xi) + \frac{t}{12}
     \Scal_M (0) + O(t^{3/2}).
 \end{equation}
 \end{Theorem}

 \begin{Proof} For a proof see \cite{Roe}, Prop. 5.25, p. 73 together with
 the computation of $u_0$ and $u_1$ on page 78. Note that we
 consider $\exp (t\Delta_M /2 )$ instead of $\exp (t\Delta_M)$.
 \end{Proof}

 \noindent Thus we have in the same local coordinates as in
 Section 5.3
 $$
     P_t f (x) = S(t) \left((\left.\det\right.^{-1/4} g^M (\xi)
     + \frac{t}{12}\Scal_M (0))f + O(t^{3/2})\right),
 $$
where $S(t)$ is the operator in Proposition \ref{pseudo}. By the expansion (\ref{metrik}) for the metric tensor we obtain
 \begin{equation}\label{g^viertel}
 \left.\det\right.^{-1/4} g^M(\xi) = 1 + \frac{1}{12}\Ric^M_{uv}(0) \phi^u (\xi)
 \phi^v (\xi) + O(\vert \xi \vert^3),
 \end{equation}
 where $\Ric^M$ denotes the {\em Ricci-tensor} of $M$. Therefore
 \begin{eqnarray}\nonumber
     \Delta \left.\det\right.^{-1/4} g^M (0) &=& \sum_{i=1}^l\frac{
     \Ric^M_{uv}}{6} \frac{\partial\phi^u}{\partial \xi^i}\frac{\partial\phi^v}{
     \partial \xi^i}(0)= \sum_{i=1}^l\langle e_i,\frac{\Ric^M}{6} e_i\rangle(x)\\
 &=&\label{Deltag^viertel} \frac{1}{6}\overline{\Ric}_{M/L}(x),
 \end{eqnarray}
 where $\overline{\Ric}_{M/L}$ denotes the partial trace of the Ricci-tensor
 using an orthonormal base $e_i$ of $\phi_* T_xL$. Evaluating (\ref{g^viertel}) 
 at $\xi = 0$ we obtain 
 \begin{eqnarray*}
     P_t f (x) &=&  e^{t\left(-\frac{
 \Scal_L}{4}+\frac{\vert\tau\vert^2}{8}+\frac{\overline{R}_{M/L}}{12}\right)(x)}(
1 +\frac{t}{12}\Scal_M (x)) f(x)  \\ 
& & + \left.
  \frac{t}{2}\Delta_L (\left.\det\right.^{-1/4} g^M f) (x) + O(t^{3/2})\right) \\
  &=& e^{-\frac{t
 \Scal_L}{4}+\frac{t\vert\tau\vert^2}{8}+\frac{t(\overline{R}_{M/L}+\Scal_M)}{12}}f(x)\\ 
& & + \frac{t}{2}( \Delta_L  +  \Delta_L \left.\det\right.^{-1/4} g^M(x))f(x) + O(t^{3/2}).
 \end{eqnarray*}
 and with (\ref{Deltag^viertel}) we finally arrive at the
 asymptotic formula
 \begin{eqnarray}\label{heatpinexp}
     & & P_t f(x)\\\nonumber  &=&  e^{t\big(-\frac{
     \Scal_L}{4}+\frac{\vert\tau\vert^2}{8}+\frac{\overline{R}_{M/L}+
     \overline{\Ric}_{M/L} +\Scal_M}{12}\big)}f(x) + \frac{t}{2}\Delta_L f(x) +
     O(t^{3/2}).
 \end{eqnarray}
 Therefore the analogues of Theorem \ref{T1} and the Remark at the end of the previous section consist of the following
 statement:

 \begin{Theorem}\label{heatpin} Using the notations from above the operator 
families $(B(t))$ and $(T(t))$ defined on $C(L)$ by 
 $$ B(t) f (x) :=\frac{e^{t\big(\frac{\Scal_L}{4}
- \frac{\vert\tau_{\phi}\vert^2}{8}-\frac{\overline{R}_{M/L}
 + \overline{\Ric}_{M/L}+\Scal_M}{12}\big)(x)}}{\sqrt{2 \pi t}^{l-m}}\int_L
 p^M_t(x,y) f(y)\vol_L(dy)
 $$
 and 
$$ T(t)f (x) = \frac{\int_L p^M_t(x,y) f(y)\vol_L(dy)}{\int_L
 p^M_t(x,y)\vol_L(dy)} $$
are both Chernoff equivalent to the heat semigroup on $L$.
 \end{Theorem}
 
\section{Limit Densities for the Pinning Construction under global 
normalization}\label{Limitdens}

 \subsection{Introduction}\label{Limitintro}

 We return to the introduction of section \ref{properex}. The operator families
 $S = (S(t))$ defined by (\ref{V}) are not proper in general. In our examples 
 their rescaled versions $T = (T(t))$ and $B = (B(t))$, cf. (\ref{U1}) and 
(\ref{U2})), become proper and Chernoff equivalent to the heat semigroup and 
 thus, applying the pinning construction with either $T$ or $B$ yields 
 (assuming an appropriate choice of the interpolation family
 $\mathcal{Q}$) weak convergence of the measures to the law $\W^x_{L}$ of 
 Brownian motion on $L$.\\
 
On the other hand we can perform the pinning construction directly with 
$S$. Then for each partition $\PP$ of $[0,1]$ the comparison of the induced 
measures on $L^{\PP}$ corresponding to $B$ and $S$ yields
 $$
     \frac{d \P^x_{L,\PP,S}}{d \P^x_{L,\PP,B}}(x; y_1,...,y_r)
      = e^{(t_1 - t_0) D(x) + \sum_{k=1}^{r-1} (t_{k+1} - t_k) D(y_k)}.
 $$
 where we use the symbol of the family as an additional index.
 If we choose the same interpolation family for both $B$ and $S$ this implies
 $$
     \frac{d \P^x_{\PP,S}}{d \P^x_{\PP,B}}(\omega)
      = e^{ (t_1 - t_0) D(x) + \sum_{k=1}^{r-1} (t_{k+1} - t_k) 
        D(\omega(t_k))}.
 $$
 on the path space $C_M[0,1]$. In the exponent, we have Riemann sums which 
 converge uniformly on compact subsets of $C_M[0,1]$ to the corresponding 
 integrals. According to Theorem \ref{mainpinning} for $|\PP_k| \rightarrow 0$
 the sequence $(\P^x_{\PP_k,B})$ is uniformly tight.
 Therefore the measures $\P^x_{\PP_k,S}$ converge weakly over the
 space $C_M[0,1]$ to a measure $\nu^x_{L}$ which is concentrated on
 $C_L[0,1]$ with density
 \begin{equation}\label{pindensity}
   \frac{d \nu^x_{L}}{d \W^x_L}(\omega) 
   = \lim_{\vert\mathcal{P}\vert\to 0}
\frac{d \P^x_{\PP,S}}{d \P^x_{\PP,B}}(\omega)  
= \exp\left(\int_0^1 D(\om(s)) ds\right).
 \end{equation}

  As a first application we get for $M = L$ a result which for geodesic interpolation on $L$ was essentially contained with a completely different proof first in \cite{Andersson-Driver}.

 \begin{Corollary}\label{AnderDriv} Let in the construction above $q$ be the 
pseudo-Gaussian kernel on $L$ 
$$q_t (x,y) := \frac{1}{(2\pi t)^{l/2}}
 e^{-\frac{d_L (x,y)^2}{2t}}. $$
 Let $\mathcal{Q}$ be an arbitrary
 interpolating family such that the sequence $(\P^x_{\mathcal{P}_k})$ is tight 
over $C_L[0,1]$. Then the limit measure $\nu_{L}^x$ corresponding to these 
kernels under the pinning construction is equivalent to Wiener measure 
$\W_L^x$ on $L$ with density
 $$
  \frac{d \nu_{L}^x}{d \W^x_{L}}(\om) = \exp\left(\frac{1}{6}\int_0^1 \Scal (\om(s))
     ds\right).
 $$
 \end{Corollary}

 We can state two more density results in
 the spirit of Corollary \ref{AnderDriv}. Namely, if we use the
 following pseudo-Gaussian density
 \begin{equation}\label{PSG}
 q_t (x,y) := \frac{1}{(2\pi t)^{l/2}} e^{-\frac{d_M (x,y)^2}{2t}}
 \end{equation}
 but integrate with respect to the volume form on $L$ we obtain
 the following result from Theorem \ref{T1} in the same way as in
 the preceding proof.

 \begin{Corollary}\label{T2}
 The family of measures obtained by the pinning construction to the
 rescaled restriction  (\ref{PSG}) of the pseudo-Gaussian kernel on $M$ 
 integrated with respect to $L$
 converges weakly to a measure $\nu_{L}^x$ which is equivalent to
 Wiener measure on $L$ with density
 \begin{eqnarray*}
     \frac{d \nu_L^x}{d \W^x_{L}}(\om) = \exp\left(\int_0^1 \left\lbrace
     \frac{1}{4}\Scal_L -
     \frac{1}{8}\vert\tau\vert^2 - \frac{1}{12}\overline{R}_{M/L}\right\rbrace
     (\om(s))
     ds\right).
 \end{eqnarray*}
 \end{Corollary}
 If we use the properly normalized heat kernel on $M$, namely
 \begin{equation}\label{PHK}
     q_t (x,y) := (2\pi t)^{(m-l)/2} p_t^M (x,y)
 \end{equation}
 we get from Theorem \ref{heatpin}:

 \begin{Corollary}\label{T3}
 The family of measures obtained by the pinning construction applied to the 
 restriction of the rescaled heat kernel (\ref{PHK}) of $M$ to $L$ converges 
weakly to a measure $\nu_L^x$ which is equivalent to Wiener measure on $L$ 
with density
 $$\frac{d \nu_L^x}{d \W^x_{L}}(\om) = \exp\left(\int_0^1 D(\om(s)) ds\right)
 $$
 where the function $D$ on $L$ is defined by
 \begin{equation}\label{heatdens}D(y) = \left(\frac{1}{4}
 \Scal_L - \frac{1}{8} \vert\tau_{\phi}\vert^2 -
 \frac{1}{12}(\overline{R}_{M/L} + \overline{\Ric}_{M/L}+\Scal_M)\right)(y).
 \end{equation}
\end{Corollary}

Note that the function $D$ in (\ref{heatdens}) is precisely the function in 
the exponent in (\ref{Dichte}) in Theorem \ref{surfacemeasure}. 
The explanation will be given below. 

 \subsection{The Pinning Construction as Conditional
 Probability}\label{conditional}

 We now want to apply our convergence results to the comparison of the Wiener
 measures on the manifolds $L$ and $M$ respectively. To do so, we compare the
 pinning construction applied to the family $S$ with another one obtained by
 normalizing at all partition times simultaneously.
 To be precise we construct from $\P^x_{L,\PP,S}$ the probability measure
 $$\label{defpps}
 \P_{L,\PP}^{x,\Sigma}(dy_1,...,dy_r):= C_{\PP}(x) \P^x_{L,\PP,S}
 (dy_1,...,dy_r),
 $$
 where the constant is chosen so that this becomes a probability
 distribution. In the special case where, similarly as in Corollary \ref{T3},
 the kernel is given by the restriction of the heat 
kernel $(p_t^M)$
 on $M$, the resulting probability measure $\P_{L,\PP}^{x,\Sigma}$
 on $L^\PP$ describes the {\em marginal distribution of the
 $M$-Brownian motion $W^x_M$ under the condition that it visits $L$
 at all times $t_i$ in the partition}. (Note that now the rescaling constants 
 in (\ref{PHK}) are no longer necessary because they are incorporated in the 
 new constant $C_{\PP}(x)$.) Since the volume measure of
 $L$ can be obtained from the volume measure on $M$ by the usual
 surface measure limiting procedure, we can compute the measure
 $\P_{L,\PP}^{x,\Sigma}$ also by
 \begin{eqnarray}\label{rdimsurf}
 &&\P_{L,\PP}^{x,\Sigma}(A_1 \times \cdots \times A_r) \\
\nonumber &=&\lim_{\eps\to 0} \frac {\int_{
 U_\eps(A_1)\times \cdots\times U_\eps(A_r) }\, p_{t_1}^M (x,
 dy_1)...p_{1-t_{r-1}}^M(y_{r-1},dy_r)} {\int_{ U_\eps(L)^r }\,
 p_{t_1}^M (x, dy_1)...p_{1-t_{r-1}}^M(y_{r-1},dy_r)}.
 \end{eqnarray}
 Thus we can consider $\P_{\PP}^{x,\Sigma}$ as a
 kind of surface measure induced by the heat kernel of $M$
 on $L^{\PP}$. Hence the upper index '$\Sigma$'. To point
 out the difference to the pinning distribution marginal we write
 the latter as
 $$
 \P^x_{L,\PP} (A_1 \times \cdots A_r)=\lim_{\eps\to 0} \frac {\int_{
 U_\eps(A_1)\times \cdots\times U_\eps(A_r) }\, p_{t_1}^M (x,
 dy_1)...p_{1-t_{r-1}}^M(y_{r-1},dy_r)} {\prod_{k=1}^{r}\int_{
 U_\eps(L)}\, p_{t_k-t_{k-1}}^M (y_{k-1}, dz)}
 $$
 letting $y_0 = x$, the difference being that in this last expression the kernel is
 renormalized at each partition time.

 Let $\W^{x,\Sigma}_{\PP}$ denote the law on $C_M[0,1]$ which we get from
 $\P^{x,\Sigma}_{L,\PP}$ by using {\em $M$-Brownian bridge measures} for the
 interpolating family $\mathcal{Q}$. Then $\W^{x,\Sigma}_{\PP}$
 can be considered as the law of {\em $M$-Brownian motion} conditioned to be on the
 submanifold $L$ at the times $t_k\in\PP$. Since the normalizing map $\nu \mapsto
 \frac{1}{\nu(C_M[0,1])}\nu$ is continuous with respect to the topology of 'weak convergence' on the cone of finite positive measures and since the convergence
 of the nonnormalized measures is known from Corollary \ref{T3} we get finally 
the following reformulation of Theorem \ref{surfacemeasure}: Let $x\in L$. 
 As the mesh $\vert\PP\vert$ of a partition $\PP$ of $[0,1]$ converges to $0$,
 the conditional law   
 $$
     \W^x_M (d\om \, \vert \, \om(t_k)\in L\, , \, t_k \in\PP )
 $$
 of Brownian motion on $M$, conditioned to visit $L$ at all partition times, 
 tends weakly over $C_M[0,1]$ to the measure $\mu_L^x$ which is equivalent to 
 $\W^x_L$ with density (\ref{Dichte}).

\section{Conclusion} 

Let $L$ be a smooth closed Riemannian manifold. For various one-parameter
families $S = (S(t))$ of kernel operators of the form 

\begin{equation}\label{S(t)} 
S(t)f(x) = \int_L q(t,x,y) f(y) vol_L(dy)
\end{equation}

we can verify the short time asymptotics 

\begin{equation}\label{shorttime} 
S(t) f(x) = e^{tD(x)} + \frac{t}{2} \Delta_Lf(x) + O(t^{3/2})
\end{equation}

for all $f \in C^3(L)$ and some function $D \in C(L)$ which depends on $S$
and is a combination of curvature terms. 

For each partition $\PP$ of $[0,1]$ and interpolation by geodesics or Brownian 
bridges the family $S$ induces a measure on the path space $C_L[0,1]$. As the 
partition gets finer we prove convergence in law of these measures 
to a measure $\nu_L^x$ which is equivalent to the law $\W^x_L$ of $L$-valued 
Brownian motion with Radon-Nikodym density 

\begin{equation}\label{RN} 
\frac{d \nu_L^x}{d \W^x_{L}}(\om) = \exp\left(\int_0^1 D(\om(s)) ds\right).
 \end{equation}

If $q(t,x,y)$ is replaced by $\tilde{q}(t,x,y) = 
\frac{q(t,s,y)}{\int_L q(t,x,y) vol_L(dy)}$ then the associated probability 
measures on the path space converge to $\W^x_L$. 

If $q(t,x,y)$ is the restriction of the heat kernel on a surrounding manifold 
$M$ to $L$ then the above results imply the convergence of the conditional 
law of Brownian motion on $M$, conditioned to return to $L$ at all partition 
times to a measure of the form (\ref{RN}) with an explicitly known funtion 
$D$.  Here $\nu_L^x$ can be viewed as the infinite dimensional surface measure
induced by $\W_M^x$ on the set $C_L[0,1] \subset C_M[0,1]$. 

The key tools for these results are a) asymptotic computation of Gaussian 
integrals combined with b) various differential geometric computations for 
the proof of  (\ref{shorttime}) and c) a new version of Chernoff's 
theorem in semigroup theory combined with d) tightness results for the 
convergence of the measures.

\end{article}
\end{document}